
\documentstyle{amsppt}
\document
\magnification=\magstep1
\pagewidth{6.5 true in}
\pageheight{9 true in}
\NoBlackBoxes
\define \blank{{\vskip 10pt}}
\define \dash{\text{--}}
\define \Hom{\operatorname{Hom}}
\define \pp{\bold P^2} 
\define \ppdual{\check{\bold P}^2} 
\define \pl{\bold P^1}	

\define \fpos{\underset S \to \times} 
\define \pd#1#2{\frac{\partial #1}{\partial #2}}
\define \pnd{\bold P^{N(d)}}
\define \fpopn{\underset {\pnd} \to \times}
\define \fn{F(\bold n)}
\define \cn{\Cal C_{\pnd}(\bold n)}
\define \sn{\Cal S(\bold n)}
\define \fnopen{\fn^{\circ}}
\define \cnopen{\cn^{\circ}}
\define \snopen{\sn^{\circ}}
\define \fnplus{\fn^+}
\define \cnplus{\cn^+}
\define \fnq{\fn^q}
\define \cnq{\cn^q}
%
\define \liftcusp{1} 
\define \seq{2}
\define \pandq{3} 
\define \quadrln{4} 
\define \newbasis{5} 
\define \maindiagram{6}
\define \eyesphis{7} 
\define \focalcherns{8} 
\define \infinityquad{9} 
\define \indet{10} 
\define \pcn{11} 
\define \indcnf{12} 
\define \ufeqn{13} 
\define \matpart{14} 
\define \col{15}
\define \base{16}
\define \generators{17}
\define \aufeqn{18} 
\define \matparttwo{19}
\define\gp{20} 
\define \nonsing{21} 
\define \triptwo{22}
\define \dimineq{23} 
%
%
%
\define \aone{1} 
\define \atwo{2} 
\define \athree{3} 
\define \ass{4} 
\define \ckho{5} 
\define \cktrip{6} 
\define \collino{7} 
\define \elb{8} 
\define \esone{9} 
\define \estwo{10} 
\define \fkm{11} 
\define \hilbert{12} 
\define \dj{13} 
\define \ktrans{14} 
\define \kmult{15} 
\define \ksone{16} 
\define \kstwo{17} 
\define \ksthree{18} 
\define \maillard{19} 
\define \ms{20} 
\define \mattuck{21} 
\define \mx{22} 
\define \rsone{23} 
\define \rstwo{24} 
\define \schub{25} 
\define \schubtwo{26} 
\define \semple{27} 
\define \stanley{28} 
\define \vain{29} 
\define \math{30} 
\define \zeuthen{31} 
\leftheadtext{S. J. Colley and G. Kennedy}
\rightheadtext{Simultaneous Higher-order Contact}
\topmatter
\title The enumeration of \\
simultaneous higher-order contacts \\
between plane curves
\endtitle
\author Susan Jane Colley and Gary Kennedy \endauthor
\affil Oberlin College\\
Ohio State University at Mansfield
\endaffil
\address Department of Mathematics, Oberlin
College, Oberlin,  Ohio 44074, U.~S.~A.\endaddress
\address Ohio State University at Mansfield, 1680 University Drive,
Mansfield, Ohio 44906, U.~S.~A. \endaddress
\abstract
Using the Semple bundle construction of [\ckho] and [\cktrip], we derive an
intersection-theoretic formula for the number of simultaneous contacts of
specified orders between members of a generic family of degree $d$ plane
curves and finitely many fixed curves.  The contacts counted by the formula
occur at nonsingular points of both the members of the family and the fixed
curves.
\endabstract
\thanks
During the preparation of this paper, the second author
was also an Affiliate Scholar of Oberlin College.  We wish to thank
Steve Kleiman and Bob Speiser for their suggestions and encouragement.
We are also grateful for the referee's thoughtful and very helpful comments.
\endthanks
\endtopmatter
\heading Introduction
\endheading
\smallpagebreak
\par
Two plane algebraic curves are said to have
{\it contact  of order} $o$ at  a common point  $P$
if each curve is  smooth at
$P$  and if the intersection  number at  $P$  is   $o$.
Thus, for example, a contact of order 1
is a transverse intersection, a contact of order 2
(sometimes called an {\it ordinary contact})
is a point of tangency;
more generally,
a contact of order $o$ is a point at which,
in appropriate local coordinates,
the Taylor expansions of the curves
agree up to order $o-1$.
\par
In [\fkm], Fulton, Kleiman, and MacPherson consider,
{\it inter alia},
a $p$-parameter family of plane curves together with
$p$ individual curves.
They compute, in terms of certain ``characteristic numbers'',
the number of members of the family which
simultaneously have an ordinary contact with each of the
curves.  The analysis has two parts, the first of which is
a formal calculation
of an intersection number.
The second part consists in establishing that,
under stipulated hypotheses,
this intersection number and the characteristic numbers
have their intended geometric meanings.
(The tendentious description we have just given greatly understates
the scope of the contact formula of [\fkm].
The contact formula for plane curves is an implicit special case---see
further remarks in the next paragraph.)
\par
In the present paper, by a similar two-part procedure, we
derive a higher-order contact formula.
We suppose that  $C_1, C_2,\dots, C_p$ are plane curves,
and that
$\Cal X$ is an $s$-parameter family of plane curves of fixed degree.
We suppose that  $o_1, o_2,\dots, o_p$
are specified positive integers, with sum equal to $s+p$.
Our formula then counts, under stipulated hypotheses,
the number of members of the family
which simultaneously have a contact of order
$o_1$ with $C_1$,  a contact of order $o_2$ with $C_2, \dots,$
and a contact of order $o_p$ with $C_p$.
The formula is both a specialization
and a generalization of the general contact formula of [\fkm]:
a specialization in that we consider only plane curves,
a generalization in that it counts contacts of arbitrary order.
Our formula, like that of [\fkm],
is stated using the formalism of ``modules'',
a notion which we explain in \S 3.
\par
The principal tool used in the proof, and in defining appropriate
higher-order characteristic numbers, is a tower of $\pl$-bundles
implicitly introduced by Semple [\semple] and later reintroduced
by Collino [\collino].  These bundles, which parametrize in a particularly
lucid way the higher-order curvilinear data of the plane,
are naturally adapted to the study of higher-order contact.
They generalize the variety of second-order data of $\pp$
studied by J.~Roberts and R.~Speiser in [\rsone] and [\rstwo], and, in turn,
are currently being generalized by E.~Arrondo, I.~Sols, and R.~Speiser
[\ass].  (They are developing a theory of ``derived triangles" that
provides a general approach for obtaining the higher-order data,
both curvilinear and higher-dimensional, of any scheme.)
The Semple bundles conform to the intuitive heuristic that
$(n+1)$st-order curvilinear data should fiber over the $n$th-order
data.  Moreover, since the $(n+1)$st Semple bundle variety
is a $\pl$-bundle over the
$n$th variety, the intersection rings are straightforward to calculate.
Thus the Semple bundles are superior to the canonical parameter
space for data, namely the Hilbert scheme $\operatorname{Hilb}_n\pp$
of zero@-dimensional subschemes of $\pp$,
in this respect.  Indeed, considerable
effort has been devoted recently
to determining the homology and cohomology
of punctual Hilbert schemes and to using these results to develop
enumerative applications.
(See, for example, [\elb], [\esone], [\estwo], [\ms].)
We have not explored the connections
between Semple bundles and punctual
Hilbert schemes in any detail, however.
Nor have we begun to understand the
apparently close connection between Semple bundles and
certain subschemes
of Kleiman's iterative multiple-point schemes [\kmult].
\par
In the first two sections of this paper we describe
the Semple bundle varieties
and their intersection rings.  In \S 3 we obtain a
``proto-contact formula'', and in \S 4 we show that under
stipulated hypotheses this formula counts, as intended,
the number of simultaneous contacts between a generic family
and specified curves.    The proof involves a detailed analysis
of the relationship between the universal family of
plane curves of degree $d$ and the Semple bundles which should
be of independent interest.  The fifth section briefly discusses the
ingredients of the contact formula, i\.e\., the ``higher-order
characteristic numbers''.  The final section treats special cases
and variants of the contact formula, compares the formula with
those of de~Jonqui\`eres and Fulton-Kleiman-MacPherson,
and presents an example.
\par
We assume that all our varieties and schemes are defined over an
algebraically closed
field of characteristic zero.
However, the results presented remain true with only
minor modifications
if the characteristic is sufficiently large, provided
one recognizes that in
positive characteristic our intersection numbers only
count weighted numbers of contacts.
\par
\bigpagebreak
\heading 1. Semple bundles
\endheading
\smallpagebreak
\par
We briefly recall here the definition of Semple's bundles
of higher-order curvilinear data;
for further discussion see [\ckho] or [\cktrip].
The inductive construction begins by declaring that
$F(0)$ is the projective plane and that
the first {\it Semple bundle variety}
$F(1)$ is ${\bold P}TF(0)$,
the total space of the projectivized tangent bundle.
(We use ${\bold P}E$ for the variety
representing rank 1 subbundles of the vector bundle $E$,
rather than for the variety representing rank 1 quotient bundles.)
Inductively suppose that
$F(n)$ is the projectivization of a rank 2
subbundle of the tangent bundle $TF(n-1)$.
The  {\it focal plane} at a point $p \in F(n)$
is  defined to be the
preimage,   via
the   derivative  of   the  projection $f_n : F(n) \to F(n-1)$,
of the line in  $T_{f_n (p)}F(n-1)$ represented by $p$.
This  construction gives
rise  to a rank 2 {\it bundle of focal  planes}, denoted
$\Cal F_n$;
we define $F(n+1)$
to be the total space of the projectivization
of this bundle.
In this manner we obtain a tower of $\pl$-bundles:
\smallpagebreak
\par
$$
\minCDarrowwidth{0pt}
\CD
:\\
@VVf_4V \\
F(3)   \\
@VVf_3V \\
F(2)   \\
@VVf_2V \\
F(1)   \\
@VVf_1V \\
F(0)
\endCD
$$
\smallpagebreak
\par
Suppose that $C$ is a reduced plane curve,
and that $p \in C$ is a nonsingular point.
The fiber of $F(1)$ over $p$
parametrizes the various tangent directions at $p$,
including the tangent direction $p_1$ of $C$.
The totality of such tangent directions to $C$
form a curve in $F(1)$;
the closure of this curve is called the
{\it first lift} of $C$, and denoted $C(1)$.
Now observe that the tangent direction $p_2$
of $C(1)$ at $p_1$ maps,
via the derivative of the projection
$f_1:F(1) \to F(0)$,
to $p_1$, the tangent direction of $C$ at $p$.
Hence the {\it second lift} $C(2)$,
i\.e\., the lift of the lift of $C$,
is in fact a curve in $F(2)$.
The same argument shows that there are
{\it higher-order lifts}
$C(3) \subset F(3)$,
$C(4) \subset F(4)$, etc.
We call the rational map
$\lambda$ from $C$ to $F(n)$ the {\it lifting map};
it is regular away from the singularities of $C$
and a birational map to $C(n)$.
We note that the $n$th lift $C(n)$ is just the $n$th
blow-up of $C$ at its singular points.  However,
for our enumerative purposes we need to understand
the embedding of $C(n)$ in the Semple bundle variety
$F(n)$.
\par
It is sometimes convenient to regard a point of $F(n)$
as an equivalence class of irreducible germs of plane
curves.  Thus we sometimes say that the irreducible germ of
a curve $C$ {\it represents} the point $p_n \in F(n)$,
meaning that the lift of  $C$ passes through $p_n$
above the closed point of the germ.
\par
To illustrate the transparent nature of calculations in these
Semple bundles, let us consider
a reduced plane curve $C$ defined,
in an affine chart with
coordinates $x$ and $y$,
by $f(x,y)=0$.
Then one can show that on $F(n)$
there is a {\it primary chart} isomorphic to affine $(n+2)$-space,
with
coordinates $x, y, y', y'', y^{(3)},\dots, y^{(n)}$,
and that
the ideal defining
the $n$th lift $C(n)$
includes the sequence of functions
$f, f', f'', f^{(3)},\dots, f^{(n)}$
obtained by
repeated implicit differentiation
with respect to $x$.
(There is a second primary chart in which the roles
of $x$ and $y$ are reversed.)
\par
In general these functions may not
generate the ideal defining the lift.
Consider, for example, the cuspidal cubic
$y^2=x^3$.
Then the first lift is defined in the primary chart
by
$$
y^2=x^3, \quad
2yy'=3x^2, \quad
x = \frac{4}{9}(y')^2, \quad
\text{and} \quad
y=\frac{8}{27}(y')^3. \tag\liftcusp
$$
The variety defined by just the first two equations
contains a spurious component over the origin.
(To obtain the third equation, we square both sides of
the second equation, then use the first equation to
replace $y^2$; since the curve is singular at the origin
it is then legal to cancel $x^3$ from both sides.  The
fourth equation is obtained in a similar fashion.
Clearly the third and fourth equations define an
irreducible curve in $\bold A^3$,
and the first two equations are redundant.)
\par
If we continue this example one step further we see
that the lift of a curve may leave the primary charts.
Indeed, if we implicitly differentiate the third equation
of \thetag{\liftcusp},
we obtain this equation for $C(2)$:
$$
1=\frac{8}{9}y'y''.
$$
Since the unique point $p_1$ of $C(1)$
over the origin is $(x,y,y')=(0,0,0)$,
and since $C(2)$ must be complete,
the unique point $p_2$ of $C(2)$ over $p_1$
must be the one point
over $p_1$ missed by the primary chart,
the point which, intuitively, represents
infinite curvature.
There is one such point on each fiber
of $F(2)$ over $F(1)$;
taken together, these points form a section of
the $\pl$-bundle
which we call the {\it divisor at infinity} and denote by $I_2$.
\par
Here is another characterization of the divisor at infinity.
The derivative of $f_1$ at $p_1$,
$$
df_1:T_{p_1}F(1) \to T_{f_1(p_1)} F(0),
$$
maps a 3@-dimensional vector space to a 2@-dimensional vector space.
The kernel is one-dimensional, so there is a unique
direction annihilated by $df_1$.
This direction is represented by a point on the divisor at infinity.
More generally, the derivative of $f_n$ at a point $p_n$,
$$
df_n:T_{p_n}F(n) \to T_{f_n(p_n)} F(n-1),
$$
likewise has a one-dimensional kernel.
Hence there is a divisor at infinity $I_{n+1}$ on $F(n+1)$;
intuitively, a point on this divisor represents (in addition to
certain lower-order data)
infinite curvilinear data of order $n+1$.
To avoid clumsy notation, we will also denote by
$I_{n+1}$ the pullback of this divisor to any  Semple bundle variety
above $F(n+1)$ in the tower, and continue to call the pullback a divisor
at infinity.
Note that the divisors at infinity have normal crossings; in particular,
the intersection of two such divisors has codimension two.
\par
The simplest sort of chart meeting a divisor at infinity is
a {\it secondary chart}.
To specify such a chart on $F(n)$
-- over the affine chart of $\pp$
with coordinates $x$ and $y$ --
choose an integer
$j$  between 2 and $n$.
Then there is a chart isomorphic to affine $(n+2)$-space,
with coordinates
$x, y, y', y'',\dots,y^{(j-1)}, x', x'',\dots,x^{(n-j+1)}$.
For $i=1,\dots,j-1$, the coordinate $y^{(i)}$ measures
$dy^{(i-1)}/dx$, a ratio of differentials of coordinates on $F(i-1)$;
for $i=1,\dots, n-j+1$, the coordinate $x^{(i)}$ measures
$dx^{(i-1)}/dy^{(j-1)}$, a ratio of differentials of coordinates
on $F(j+i-1)$.
(See [\ckho] for an explanation of
a full system of charts.)
If $C$ is
a reduced plane curve defined
by $f(x,y)=0$, then
the ideal defining the
$n$th lift $C(n)$
includes a sequence of functions
obtained by
repeated implicit differentiation.
In the first $j-1$ of these differentiations
we treat $x$ as the independent variable,
and denote the derivative $dy/dx$ of $y$ by $y'$,
the derivative $dy'/dx$ of $y'$ by $y''$, etc.
In the remaining $n-j+1$ differentiations
we treat $y^{(j-1)}$ as the independent variable,
and denote the derivative $dx/dy^{(j-1)}$ of $x$ by $x'$,
the derivative $dx'/dy^{(j-1)}$  of $x'$ by $x''$, etc.
The derivative $dy^{(i)}/dy^{(j-1)}$ of $y^{(i)}$ ($0\leq  i < j-1$)
is obtained by the chain rule:
$$
\frac{d y^{(i)}}{d y^{(j-1)}} = \frac{d y^{(i)}}{dx}\frac{dx}{d y^{(j-1)}}
= y^{(i+1)}x'.
$$
In other words, the defining ideal
for $C(n)$ includes the functions
$$
f, P(f), P^2(f), \dots, P^{j-1}(f),
QP^{j-1}(f), Q^2P^{j-1}(f), \dots, Q^{n-j+1}P^{j-1}(f),
\tag\seq
$$
where $P$ and $Q$
are the differential operators
$$
\aligned
P &= \pd{}{x} + y'\pd{}{y} + y''\pd{}{y'}
	+ \cdots + y^{(j-1)}\pd{}{y^{(j-2)}}\\
Q &= x'P + \pd{}{y^{(j-1)}} + x''\pd{}{x'} + x^{(3)}\pd{}{x''}
	+ \cdots + x^{(n-j+1)}\pd{}{x^{(n-j)}}.
\endaligned
\tag\pandq
$$
Away from the singularities of $C$, the
sequence of functions in \thetag{\seq} generates the ideal,
but additional generators will be needed
to eliminate spurious components over
the singularities.
In this {\it $j$th secondary chart} the divisor at infinity
$I_j$ is defined by the vanishing of $x'$.
\par
To illustrate the rules for calculating in secondary charts,
we continue our example of the cuspidal cubic
$y^2=x^3$.  Implicitly differentiating the
third equation of \thetag{\liftcusp} with respect to $y'$,
we obtain
$$
x'=\frac{8}{9}y'.
$$
Here $x'$ is the coordinate of the secondary chart
measuring $dx/dy'$.  (In other words, $x'$
is the reciprocal of the ordinary second derivative coordinate
$dy'/dx$.)  From the equations of \thetag{\liftcusp}
one easily sees that the first lift is tangent to the fiber of
$F(1)$ over the origin.
Note that, as expected,
the second lift hits the
divisor at infinity $I_2$
over the origin.
\par
We note for future reference that, taken together,
the primary and secondary charts
cover all of the Semple bundle variety
except for intersections of two or more divisors at infinity.
\par
By definition $F(n+1)$ is a subvariety of ${\bold P}TF(n)$,
the total space of the projectivized tangent bundle of $F(n)$.
If $n \geq 2$, then ${\bold P}TI_n$, the total space of the
projectivized tangent bundle of the divisor at infinity,
is likewise a subvariety of ${\bold P}TF(n)$; its codimension
is 2.  Now one can easily verify, e\.g\., by a calculation in
local coordinates, that $F(n+1)$ and ${\bold P}TI_n$
are transverse.  Hence their intersection is a codimension
2 subvariety of $F(n+1)$ which we call the
{\it locus of tangency} to $I_n$.
A point of this locus represents a point of $I_n$
together with a tangent direction belonging to
the fiber of the focal plane at that point.
Again to avoid clumsy terminology, we continue
to speak of the ``locus of tangency'' when we ought
to say ``the pullback through the Semple bundle tower
of the locus of tangency''.
The $j$th secondary chart on $F(n)$ meets only one such locus,
namely the locus of tangency to $I_j$; this locus is defined
by the vanishing of both $x'$ and $x''$.  (But if $j=n$, there
is no such locus meeting the secondary chart at all.)
\par
Suppose now that $\Cal X$ is a family of plane curves,
and that its general member is reduced.
For  each reduced member there is a (rational) lifting map
to $F(n)$;  these maps fit together to form
a rational map $\lambda:\Cal X \dashrightarrow F(n)$.
We call the closure
of $\lambda(\Cal X)$ the {\it lift} of the family, and denote
it by $\Cal X(n)$.
By definition, $\Cal X(n)$ is the union of the graph of
$\lambda$ and a closed subvariety $\Cal S(n)$, each point of which lies
over a singular or nonreduced point
of the corresponding member of $\Cal X$.
\par
\bigpagebreak
\heading 2. Intersection rings
\endheading
\smallpagebreak
\par
Since the variety $F(n)$ is the projectivization of the
rank 2 bundle $\Cal F_{n-1}$ over $F(n-1)$, the intersection rings
may be determined inductively from standard theory.
Specifically, the Chow ring $A^* (F(n))$ is an
$A^* (F(n-1))$-algebra
generated by the tautological class
$\phi_n := c_1(\Cal O_{F(n)}(1))$,
which satisfies a single quadratic relation
$$
\phi_n^2 + c_1(\Cal F_{n-1})\phi_n + c_2(\Cal F_{n-1}) = 0.
\tag\quadrln
$$
(Note:  Here and in the sequel,
we omit pullbacks of classes when convenient and
when no confusion should result.)
The base variety $F(0)$ is $\pp$, whose
Chow ring is generated by
the hyperplane class $h$, subject to the relation $h^3=0$.
Thus $A^* (F(n))$ is generated by
$h$, $\phi_1,\dots,\phi_n$ subject to the relations
mentioned above.
\par
For the purposes of studying contact between plane curves,
however, this
description is not desirable.  Instead, we provide an equivalent
formulation using more geometric classes.
In particular, we will eliminate the
$\phi_n$'s in favor of $\check h$,
the dual hyperplane class on $\ppdual$, and
the classes $i_j$ of the divisors at infinity $I_j$.
\proclaim{Theorem 1}
For $n \ge 1$,
the Chow ring $A^* (F(n))$ is generated by
$h,\check h, i_2,\dots,i_n$ subject to the relations
$$
h^3 = 0, \quad  \check h^3 = 0, \quad h^2 - h\check h + \check h^2 = 0,
$$
and, for $k= 2,\dots,n$,
$$
i_k^2 = ((2k-1)h - (k+1)\check h - ki_2 - (k-1)i_3 -\dots-3i_{k-1})i_k.
$$
\endproclaim
\demo{Proof}
To begin, $F(1)$ is the incidence correspondence of $\pp$.  Thus
$F(1) \subset \pp \times \ppdual$ and it is well known that
$$
A^* (F(1)) \cong
\frac{\bold Z [h,\check h]}{(h^3,\check h^3,h^2 - h\check h + \check h^2)}.
$$
It follows from a change of basis calculation that
$$
\phi_1 = \check h - 2h.
\tag\newbasis
$$
\par
To finish the proof, we need to rewrite \thetag{\quadrln}
without using any
$\phi_k$'s.  To do this, we note first that the
focal plane bundle $\Cal F_k$,
$k\ge 1$, fits into the following commutative diagram
of exact sequences:
$$
\CD
0 @>>> T_{F(k)/F(k-1)} @>>> \Cal F_k
		@>>>\Cal O_{F(k)}(-1)  @>>> 0 \\
@. @VV=V @VVV @VVV @. \\
0 @>>>  T_{F(k)/F(k-1)}@>>> T_{F(k)}
		@>>> f_k^* T_{F(k-1)} @>>> 0 \\
\endCD
\tag\maindiagram
$$
(This is the dual of diagram 1 of [\cktrip].)  Let
$$
\sigma: \Cal O_{F(k+1)}(-1) \rightarrow f_{k+1}^*\Cal F_k \rightarrow
 f_{k+1}^*\Cal O_{F(k)}(-1)
$$
be the composite of the tautological map and the pullback to
$F(k+1)$ of the map in the top row of diagram \thetag{\maindiagram}.
Then $\sigma$
has zero locus equal to $PT_{F(k)/F(k-1)} = I_{k+1}$.
In addition, $\sigma$
defines a section of
 $\Hom(\Cal O_{F(k+1)}(-1),\allowmathbreak f_{k+1}^* \Cal O_{F(k)}(-1))$.
Hence
$$
\aligned
i_{k+1} &= [PT_{F(k)/F(k-1)}] \\
&=f_{k+1}^* c_1(\Cal O_{F(k)}(-1)) - c_1(\Cal O_{F(k+1)}(-1)) \\
&=\phi_{k+1} - \phi_k.
\endaligned
\tag\eyesphis
$$
Now the Euler sequence
$$
0 \rightarrow  \Cal O_{F(k)} \rightarrow
	  f_k^* \Cal F_{k-1} \otimes \Cal O_{F(k)}(1)
		\rightarrow T_{F(k)/F(k-1)} \rightarrow 0
$$
and the top sequence of \thetag{\maindiagram} together imply that,
for $k \ge 1$,
$$
\aligned
c_1(\Cal F_k) &= c_1(\Cal F_{k-1}) + \phi_k \\
c_2(\Cal F_k) &= 2c_2(\Cal F_{k-1}) + c_1(\Cal F_{k-1})\phi_k.
\endaligned
\tag\focalcherns
$$
Using \thetag{\eyesphis} to substitute
for $\phi_{k+1}$ in \thetag{\quadrln}
and simplifying, one finds
$$
i_{k+1}(i_{k+1} + 3\phi_k + c_1(\Cal F_{k-1})) = 0.
\tag\infinityquad
$$
\par
Note that the inductive formula for
$c_1(\Cal F_k)$ in \thetag{\focalcherns}
yields the closed-form formula
$$
c_1(\Cal F_k) = 3h + \phi_1 + \phi_2 + \dots +  \phi_k.
$$
Hence \thetag{\infinityquad} can be rewritten as
$$
i_{k+1}(i_{k+1} + 3h  + \phi_1 + \dots + \phi_{k-1} + 3\phi_k) = 0.
$$
In view of \thetag{\newbasis} and \thetag{\eyesphis},
this is equivalent to, for
$k \ge 1$,
$$
i_{k+1}^2 = ((2k+1)h - (k+2)\check h - (k+1)i_2 - ki_3 -
\dots - 3i_k)i_{k+1}. \quad \qed
$$
\enddemo
\par
In [\cktrip] we studied the action of
the projective general linear group $PGL(2)$
on the Semple bundles induced from the $PGL(2)$-action on $\pp$.
Since this action preserves incidence
in $\pp$, there is a special orbit $Z_k$
on each $F(k)$, $k \ge 1$, that is isomorphic to $F(1)$
and is represented by the
germ of a line.  (The intuition is that $Z_k$
measures ``zero data'' of orders 2
through $k$.)  Since $Z_k$ is represented
by the germ of a smooth curve,
it must
be disjoint from the divisors at infinity
$I_j$ for all $j \le k$.
Thus if $z_k := [Z_k]$, we have
$$
i_j \cdot z_k = 0, \quad j \le k.
$$
Also each $Z_k$ is a section of the $\pl$-bundle obtained by
restricting the Semple bundle to $Z_{k-1}$.  Hence the
intersection of $Z_k$ with the fiber of $F(k)$
over a point of $F(1)$ has degree 1.
Thus, for $k \ge 1$,
$$
\int_{F(k)}h^2\check h z_k = \int_{F(k)}h\check h ^2 z_k = 1.
$$
In a similar manner, it follows that
$$
\int_{F(k)}h^2\check h z_{j-1}i_j \cdots i_k =
	 \int_{F(k)}h\check h ^2 z_{j-1}i_j \cdots i_k = 1,
	\quad 2 \le j \le k.
$$
(Note that $Z_1$ is all of $F(1)$, so that $z_1=1$; in this case
the claims above are vacuous or obvious.)
\par
In view of these remarks and the calculation of
the Chow ring of $F(n)$ above, it is
not difficult to deduce the following matrix for
the intersection pairing of
$A^1(F(n))$ with $A^{n+1}(F(n))$,
in which $f_0 = f_1 = 1$ and $f_j$ denotes the
$j$th Fibonacci number:
$$
\matrix
	& \check h^2z_n & h^2i_2i_3\cdots i_n & h^2\check h i_3\cdots i_n
&
h^2\check h z_2 i_4\cdots i_n & h^2\check h z_3 i_5\cdots i_n
& \cdots & h^2 \check h z_{n-1}\\
	&&&&&&&\\
	h&1&0&0&0&0&\cdots&0	\\
	\check h&0&1&0&0&0&\cdots&0	\\
	i_2&0&-3&1&0&0&\cdots&0	\\
	i_3&0&5&-3&1&0&\cdots&0	\\
	i_4&0&-8&5&-3&1&\cdots&0	\\
	\vdots&\vdots&\vdots&\vdots&
		\vdots&\vdots&\ddots&\vdots\\
	i_n&0&(-1)^{n+1}f_{n+1}&(-1)^nf_n&(-1)^{n-1}f_{n-1}&
	(-1)^{n-2}f_{n-2}&\cdots&\, \, 1 \, .	\\
\endmatrix
$$
\par
\bigpagebreak
\heading 3. The proto-contact formula
\endheading
\smallpagebreak
\par
We associate to each reduced plane curve $C$ a sequence
$d$, $\check d$, $\kappa_2$, $\kappa_3$,~\dots
of {\it (higher-order) characteristic numbers},
beginning with the degree and class.
Each of the other numbers is defined by
$$
\kappa_j:=
	\int i_j  \cap\lbrack C(j) \rbrack ,
$$
i\.e\., as the intersection number of a lift of the curve
with a divisor at infinity; this intersection is defined on the
$j$th Semple bundle variety, or any Semple bundle variety
above it in the tower.
Since the lift of a curve cannot hit a divisor at infinity
over a nonsingular point,
these characteristic
numbers count certain sorts of singularities.
A point of $C$ over which $C(j)$ meets $I_j$
will be called a
{\it $j$th-order cusp},
and the
corresponding characteristic number $\kappa_j$
will be called the {\it number of $j$th-order cusps} on $C$.
For example, the singularity at
the origin of $y^2=x^{2j-1}$ is a $j$th-order cusp
and contributes 1 to $\kappa_j$.
Note that $\kappa_j$ may count with multiplicities;
for example, a curve  may have a singularity at which distinct
branches each contribute to $\kappa_j$.
One can easily show that repeated lifting will desingularize
a specified curve, hence that only
finitely many characteristic numbers
are non-zero.
\par
The {\it $n$th contact module} of a reduced plane curve $C$
is a certain element of the polynomial algebra
over the integers in the following indeterminates:
$$
\aligned
\Lambda_j,&\qquad 0\le j\\
\Pi_j,&\qquad 1\le j\\
\Gamma^k_j,&\qquad 2\le k\le j .
\endaligned
\tag\indet
$$
Specifically, it is
$$
\align
m_n(C) := d\Lambda_n &+ \check d\Pi_n \\
&+ (3\check d +\kappa_2)\Gamma^2_n \\
&+ (4\check d + 3\kappa_2 + \kappa_3)\Gamma^3_n \\
&+ (5\check d + 4\kappa_2 + 3\kappa_3 + \kappa_4)\Gamma^4_n \\
&+ \dots \\
&+ \left((n+1)\check d + n\kappa_2 + (n-1)\kappa_3 +
\dots + 3\kappa_{n-1} + \kappa_n\right)\Gamma^n_n.
\endalign
$$
In particular the 0th and 1st contact modules are
$$
m_0(C):=d\Lambda_0, \quad \text{and} \quad
m_1(C):=d\Lambda_1 + \check d\Pi_1 .
$$
If we assign weight $n$ to each of
$\Lambda_n$, $\Pi_n$, and $\Gamma^k_n$,
then the $n$th module is homogeneous of weight $n$.
\par
Now suppose that
$C_1, C_2,\dots,C_p$
are reduced plane curves.
Suppose that $\Cal X$ is a family of plane curves
over $S$, a parameter space of dimension~$s$;
suppose that the general member of the family
is reduced.
Suppose that
$n_1, n_2,\dots,n_p$ are specified positive integers,
with sum equal to $s$.
Let
$\fn$
be the product of Semple bundle varieties
$F(n_1)\times F(n_2)\times\dots\times F(n_p)$,
and let $\pi_1, \pi_2,\dots,\pi_p$
be the various projections to the factors.
The fiber product
$$
\Cal X_S(\bold n):=
\Cal X(n_1) \fpos \Cal X(n_2)\fpos
\dotsm \fpos \Cal X(n_p)
$$
is a subvariety of $\fn \times S$;
its fiber over a point $s$ of $S$ is the product of $p$ lifts
of the curve $X_s$ in the family.
Let $\sigma$ denote the projection of $\fn \times S$
onto its first factor.
\par
 $$
\CD
@. @. @. C_j(n_j) \\
@. @. @. @VVV \\
\Cal X(n_1) \fpos \Cal X(n_2)\fpos
\dotsm \fpos \Cal X(n_p)
= \Cal X_S(\bold n)
@>>> \fn \times S @>\sigma>> \fn @>\pi_j>> F(n_j) \\
\endCD
$$
\par
We define the
{\it proto-contact number}
of {\it type} $(n_1, n_2,\dots,n_p)$
by
$$
I:=\int_{\fn}
	\pi_1^*\lbrack C_1(n_1)\rbrack
	\cdot
	\pi_2^*\lbrack C_2(n_2)\rbrack
	\cdot
	\dots
	\cdot
	\pi_p^*\lbrack C_p(n_p)\rbrack
	\cap
	\sigma_*\lbrack \Cal X_S(\bold n) \rbrack .
\tag\pcn
$$
\proclaim{Theorem 2}
The proto-contact number is obtained by
multiplying the contact modules
$m_{n_1}(C_1)$, $m_{n_2}(C_2),\dots, m_{n_p}(C_p)$,
evaluating each monomial in the resulting product, and
performing the indicated arithmetic.  Each monomial in the
product is of weight $s$ and of the form
$y_1y_2\dotsm y_p$,
where each $y_j$ is either
$\Lambda_{n_j}$ or $\Pi_{n_j}$ or some $\Gamma^k_{n_j}$;
to evaluate this monomial means to replace it by
$$
\int_{\fn}
	\pi_1^*(\overline y_1)
	\cdot
	\pi_2^*(\overline y_2)
	\cdot
	\dots
	\cdot
	\pi_p^*(\overline y_p)
	\cap
	\sigma_*\lbrack \Cal X_S(\bold n) \rbrack ,
	\tag\indcnf
$$
where
$$
\overline y_j:=\left\{
\aligned
	h \qquad & \text{if } y_j=\Lambda_0,\\
	\check h^2 z_{n_j}\qquad &\text{if } y_j=\Lambda_{n_j},
		\quad n_j>0,\\
	h^2 i_2 i_3 \dotsm i_{n_j}\qquad &\text{if } y_j=\Pi_{n_j},\\
	h^2 \check h z_{k-1} i_{k+1} i_{k+2} \dotsm i_{n_j}
		\qquad &\text{if } y_j=\Gamma^k_{n_j}.
\endaligned
\right.
$$
\endproclaim
\demo{Proof}
Our discussion in \S 2 shows that the set
$$
\align
\{
&{\check h}^2 z_n,
h^2 i_2 i_3 \dotsm i_n, \\
&h^2 \check h i_3 i_4 \dotsm i_n,
h^2 \check h z_2 i_4 i_5 \dotsm i_n,
h^2 \check h z_3 i_5 i_6 \dotsm i_n, \\
&\quad
\dots,
h^2 \check h z_{n-1}
\}
\endalign
$$
forms a basis for $A^{n+1}(F(n))$,
and that the dual basis for $A^1(F(n))$ is
$$
\align
\{
&h,
\check h, \\
&3\check h +i_2,
4\check h + 3i_2 + i_3,
5\check h + 4i_2 + 3i_3 + i_4, \\
&\quad
\dots,
(n+1)\check h + ni_2 + (n-1)i_3 +\dots + 3i_{n-1} + i_n
\}.
\endalign
$$
(For $n=0$ the singleton $\{h\}$ is a self-dual basis.)
The characteristic numbers of a reduced plane curve $C$
are determined by
$$
\align
d &= \int_{\pp} h  \cap\lbrack C \rbrack
	=  \int_{F(n)} h  \cap\lbrack C(n) \rbrack,\\
\check d &= \int_{F(1)} \check h  \cap\lbrack C(1) \rbrack
	=  \int_{F(n)} \check h  \cap\lbrack C(n) \rbrack,\\
\kappa_j &= \int_{F(j)} i_j  \cap\lbrack C(j) \rbrack
	=  \int_{F(n)} i_j  \cap\lbrack C(n) \rbrack \quad (j \le n).
\endalign
$$
Hence the
rational equivalence class of the $n$th lift of $C$
is given by
$$
\align
\lbrack C(n) \rbrack =& d{\check h}^2 z_n
	+ \check d h^2 i_2 i_3 \dotsm i_n \\
&+ (3\check d +\kappa_2)h^2 \check h i_3 i_4 \dotsm i_n \\
&+ (4\check d + 3\kappa_2 + \kappa_3)
	h^2 \check h z_2 i_4 i_5 \dotsm i_n \\
&+ (5\check d + 4\kappa_2 + 3\kappa_3 + \kappa_4)
	h^2 \check h z_3 i_5 i_6 \dotsm i_n \\
&+ \dots \\
&+ \left((n+1)\check d + n\kappa_2 + (n-1)\kappa_3 +
	\dots + 3\kappa_{n-1} + \kappa_n\right)
	h^2 \check h z_{n-1}.
\endalign
$$
(For $n=0$, the class of $C(0)=C$ is of course $dh$.)
The theorem follows immediately from this formula.
\qed
\enddemo
\par
\bigpagebreak
\heading 4. The simultaneous contact formula
\endheading
\smallpagebreak
\par
Suppose that $C$ and $X$ are reduced plane curves.
A  {\it contact} (or {\it honest contact}) of {\it order} $o$
between them is a point $x \in C(o-1) \cap X(o-1)$
whose image in $\pp$ is a nonsingular point on each curve.
Note that, for nonsingular curve germs, the
following statements are equivalent:
\roster
\item"$\bullet$" There is a contact of order $o$ between them.
\item"$\bullet$" In appropriate local coordinates, the
Taylor expansions agree up to order $o-1$.
\item"$\bullet$" The intersection number is at least $o$.
\endroster
We will call a point $x \in C(o-1) \cap X(o-1)$ whose image in $\pp$
is a singular point on $C$ or $X$ a {\it false contact}.
\par
Next suppose that
$C_1,C_2,\dots,C_p$
are reduced plane curves.
A {\it simultaneous contact} of {\it order}
$(o_1,o_2,\dots,o_p)$
between $X$ and these $p$ curves
is a $p$-tuple
$(x_1,x_2,\dots,x_p)$
in which the point $x_1 \in F(o_1-1)$ is a contact of order $o_1$
between $X$ and $C_1$,
the point $x_2 \in F(o_2-1)$ is a contact of order $o_2$
between $X$ and $C_2$, etc.
\par
We say that a plane curve
$C$ has a {\it profound cusp} if, for some $j$,
the $j$th lift $C(j)$ meets the intersection of $I_j$
and the pullback of another divisor at infinity;
i\.e\., $C$ has a profound cusp at $P$
if some branch of $C$ has simultaneously a $j$th-order cusp
and a cusp of lower order.  For example, $y^3=x^5$ has a
profound cusp; its lifts meet both $I_2$ and $I_3$.
We say that
$C$ has a {\it flat cusp} if the $j$th lift of $C$
meets the locus of tangency to the
divisor at infinity $I_{j-1}$.  If $C(j-1)$ happens to be
nonsingular over $P$,
then $C$ has a flat cusp if and only if $C(j-1)$ is tangent
to $I_{j-1}$.
For example, the second lift of
$y^3 = x^4$
is tangent to $I_2$; hence this curve has
a flat cusp at the origin.
\par
A family of plane curves of degree $d$
over a parameter space $S$ determines, and is determined by,
a morphism from $S$
to the projective space $\pnd$
parametrizing such curves,
where  ${N(d)}=\binom {d+2}2 - 1$.
We call such a family {\it generic}
if it is obtained from one particular specified
family by composing  the morphism from
$S$ to $\pnd$
with a generic motion of the projective space.
In other words, ``proposition $\Cal P$ is valid for a generic family''
means that if $\Cal X$ is a family of curves over $S$
determined by
$\sigma :S\to \pnd$,
then, for all $\gamma$ in some Zariski open dense subset
of the projective linear group
$PGL({N(d)})$,
proposition $\Cal P$ is valid for the family
determined by
$\gamma \circ \sigma$.
\par
\proclaim{Theorem 3}
Suppose that $\Cal X$ is a generic family of
degree $d$ plane curves
over $S$, a parameter space of dimension
$s=o_1+o_2+\dots +o_p-p$.
Suppose that each $o_i>1$,
and that $d + 1 \geq \sum_{i=1}^p o_i$.
Suppose that  $C_1,C_2,\dots,C_p$
are reduced plane curves, none of which has a
profound cusp or a flat cusp.
Suppose that these curves have only pairwise
transverse intersections.
Then the number of simultaneous contacts
of order
$(o_1,o_2,\dots,o_p)$
between some reduced member of $\Cal X$ and
$C_1,C_2,\dots,C_p$
is the proto-contact number of type
$(o_1-1,o_2-1,\dots,o_p-1)$.
\endproclaim
\par
We have assumed that each $o_i>1$ only to avoid
a clumsy exposition.
For remarks concerning the omitted cases, see \S 6(c).
\par
To prove Theorem~3, we begin in Lemmas A and B by
analyzing a lift of the universal family
$\Cal C$ of degree $d$
plane curves over $\pnd$; a point of $\Cal C$
is an ordered pair $(P,C)$ consisting of a degree $d$ plane curve
$C$ and a point $P \in C$.  Over the affine chart with
coordinates $x$ and $y$, the hypersurface $\Cal C$ is defined
in $\bold A^2 \times \pnd$ by
$$
f(x,y) = \sum_{u+v \le d} a_{uv}x^uy^v=0. \tag\ufeqn
$$
Since we wish to study simultaneous contacts,
in Lemmas C and D we analyze the fiber product of several lifts
of the universal family.
Finally, in order to ultimately rule out the possibility
that our proto-contact formula counts false contacts,
in Lemmas E and F we analyze the higher-order data
carried by a singular point of a curve.
Our analysis uses, {\it faute de mieux}, explicit and
detailed calculations.
\par
\proclaim{Lemma A}
Suppose that $d \geq n$.
\roster
\item"$\bullet$"
Except possibly above
intersections of two or more divisors at infinity
and above the loci of tangency to divisors at infinity,
the $n$th lift $\Cal C(n)$
is smooth over $F(n)$.
\item"$\bullet$"
Over the primary chart with
coordinates $x, y, y', y'', y^{(3)},\dots, y^{(n)}$,
the $n$th lift $\Cal C(n)$ is defined
in  $\bold A^{n+2} \times \pnd$
by an ideal generated by the function $f$ of
\thetag{\ufeqn}
and the functions $f', f'', f^{(3)},\dots,f^{(n)}$
obtained by
repeated implicit differentiation with respect to $x$.
Over each point of the primary chart the matrix of
this system of linear equations has rank $n+1$.
\item"$\bullet$"
Over a point on the exceptional divisor
of the secondary chart with coordinates
$x, y, y', y'',\dots, y^{(j-1)},x', x'',\dots,x^{(n-j+1)}$
(i\.e\., a point at which $x'$ vanishes)
the $n$th lift $\Cal C(n)$ is defined
by an ideal generated by the function $f$ of
\thetag{\ufeqn}
and the sequence of functions
\thetag{\seq} obtained by the procedure explained in \S 1.
Over this point the matrix of
this system of linear equations has rank $n+1$.
\endroster
\endproclaim
\demo{Proof}
Recall that, taken together,
the primary charts cover
all of $F(n)$
except divisors at infinity,
and that
the primary and secondary charts cover
all of $F(n)$
except intersections of two or more divisors at infinity.
And note that the codimension of $\Cal C(n)$
in  $\bold A^{n+2} \times \pnd$
is $n+1$.
Hence the first claim follows from the other two.
\par
Consider the primary chart.
Let $M$ denote the $(n+1)\times({N(d)}+1)$ matrix of
partial derivatives of the functions
$f, f', f'', f^{(3)},\dots,f^{(n)}$
with respect to each of the
$a_{uv}$.
(Equivalently, since all these functions
are linear in the $a_{uv}$'s,
$M$ is the matrix of the system of linear equations.)
If the rank of $M$ at a point of $\Cal C(n)$
is $n+1$, then at this point
$f, f', f'', f^{(3)},\dots, f^{(n)}$
generate the defining ideal,
and the
projection to $F(n)$ is smooth.
Each column of $M$ begins with a monomial
in $x$ and $y$, and the other entries are  obtained
by repeated implicit differentiation:
$$
\bmatrix
\pd{f}{a_{uv}} \\
\pd{f'}{a_{uv}} \\
\vdots \\
\pd{f^{(n)}}{a_{uv}}
\endbmatrix
=
\bmatrix
x^u y^v \\
\frac{d(x^u y^v)}{dx} \\
\vdots \\
\frac{d^n(x^u y^v)}{dx^n}
\endbmatrix .
$$
To see that $M$ has rank $n+1$ at every point, consider the square
submatrix consisting of the partial derivatives of
$f, f', f'',\dots, f^{(n)}$
with respect to
$a_{00}, a_{10}, a_{20},\dots, a_{n0}$:
$$
\bmatrix
1 & x & x^2 & \cdots & x^n \\
0 & 1 & 2x  & \cdots & nx^{n-1} \\
0 & 0 & 2    & \cdots & n(n-1)x^{n-2} \\
\vdots & \vdots & \vdots & \ddots & \vdots \\
0 & 0 & 0 & \cdots & n!
\endbmatrix .
\tag\matpart
$$
\par
Now consider the secondary chart.
Let $M$ denote the $(n+1)\times({N(d)}+1)$ matrix of
partial derivatives of the functions in
\thetag{\seq}
with respect to each of the
$a_{uv}$.
(Once again we may interpret $M$
as the matrix of a system of linear equations.)
Each column of $M$ begins with a monomial
in $x$ and $y$, and the other entries are  obtained
by repeated implicit differentiation:
$$
\bmatrix
x^u y^v \\
P(x^u y^v) \\
\vdots \\
P^{j-1}(x^u y^v) \\
QP^{j-1}(x^u y^v) \\
\vdots \\
Q^{n-j+1}P^{j-1}(x^u y^v) \\
\endbmatrix .
\tag \col
$$
If the rank of $M$ at a point of $\Cal C(n)$
is $n+1$, then at this point
the functions in
\thetag{\seq}
generate the defining ideal,
and the
projection to $F(n)$ is smooth.
\par
To complete our analysis of this matrix,
we need the following formulas
concerning the polynomial ring in
the variables $x, y, y', y'',\dots, y^{(j-1)},x', x'',\dots,x^{(n-j+1)}$.
\par
\proclaim{Lemma B}
\roster
\item"(a)"
For each integer $k \geq j$,
$$
Q^i P^{j-1}(x^k) \equiv \left\{
\aligned
	&0 \quad \text{if }  i < 2(k - j + 1)\\
	&\text{a positive multiple of  } (x'')^{k - j + 1}\\
	&\phantom{0}\quad \text{if }  i = 2(k - j + 1)
\endaligned
\right.
$$
modulo the ideal generated by $x$ and $x'$.
\item"(b)"
For each pair of positive integers  $i$  and  $m$,
$$
\lbrack Q^i(x^m) \rbrack_{x=0}
= Q \lbrack Q^{i-1}(x^m) \rbrack_{x=0}
+ m x'  \lbrack Q^{i-1}(x^{m-1}) \rbrack_{x=0},
$$
where $[\;]_{x=0}$ indicates evaluation.
\item"(c)"
For each nonnegative integer $h \leq j-1$
and each nonnegative integer $i$,
$$
Q^i(y^{(h)})
=
\sum_{b=h}^{j-1}\frac{1}{(b-h)!}y^{(b)}
\lbrack Q^i(x^{b-h}) \rbrack_{x=0}
+ \text{ a polynomial in } x, x',\dots,x^{(n-j+1)}.
$$
\item"(d)"
For each nonnegative integer $h \leq j-1$,
$$
Q^i P^h(y) - \sum_{b = 0}^{j-1}\frac{1}{b!}y^{(b)}Q^i P^h(x^b)
\equiv \left\{
\aligned
	&0 \quad \text{if }  i < 2(j-1-h)+1\\
	&\text{a positive multiple of  } (x'')^{j-1-h}\\
	&\phantom{0}\quad \text{if }  i = 2(j-1-h)+1
\endaligned
\right.
$$
modulo the ideal generated by $x$ and $x'$.
\item"(e)"
For each polynomial $\phi$ in $x$ and $y$,
for each $h \geq 1$, and for each $i \geq 2$,
$$
Q^i P^h(x \cdot \phi)
\equiv h Q^i P^{h-1}(\phi)
+ \sum_{a=0}^{i-2} \binom{i}{a} Q^a P^h(\phi) \cdot x^{(i-a)}
$$
modulo the ideal generated by $x$ and $x'$.
\item"(f)"
For each nonnegative integer $k$
and each nonnegative integer $h \leq j-1$,
$$
Q^i P^h(x^k y) - \sum_{b = 0}^{j-1}\frac{1}{b!}y^{(b)}Q^i P^h(x^{k+b})
\equiv \left\{
\aligned
	&0 \quad \text{if }  i < 2(k + j - 1 - h) + 1\\
	&\text{a positive multiple of  } (x'')^{k+j-1-h}\\
	&\phantom{0}\quad\text{ if }  i = 2(k + j - 1 - h) + 1
\endaligned
\right.
$$
modulo the ideal generated by $x$ and $x'$.
\item"(g)"
For each nonnegative integer $k$,
$$
Q^i P^{j-1}(x^k y) - \sum_{b = 0}^{j-1}\frac{1}{b!}y^{(b)}Q^i P^{j-1}(x^{k+b})
\equiv \left\{
\aligned
	&0 \quad \text{if }  i < 2k + 1\\
	&\text{a positive multiple of  } (x'')^k\\
	&\phantom{0}\quad \text{if }  i = 2k + 1
\endaligned
\right.
$$
modulo the ideal generated by $x$ and $x'$.
\endroster
\endproclaim
\par
\demo{Proof}
To prove statement (a), declare the weight of $x^{(t)}$ to be $2-t$.
Then $Q$ decreases the weight of a polynomial in
$x, x',\dots,x^{(n-j+1)}$ by 1, so that
$$
Q^i P^{j-1}(x^k) = \frac{k!}{(k-j+1)!} Q^i (x^{k-j+1})
$$
is a polynomial of degree $k-j+1$ and weight $2(k - j + 1)-i$.
If $i < 2(k - j + 1)$ then the weight is positive, so each term
involves either $x$ or $x'$.  If $i = 2(k - j + 1)$ then the weight
is zero, so the only term not involving $x$ or $x'$ is a multiple
of $(x'')^{k - j + 1}$.  Clearly the coefficient of this term is positive.
\par
Observe that
$$
Q^i(x^m) = \sum x^{(e_1)} x^{(e_2)} \cdots x^{(e_m)},
$$
where the sum is over all functions  $s$ from $\{1,\dots,i\}$
to $\{1,\dots,m\}$, and $e_k$ is the number of elements of
$\{1,\dots,i\}$ for which the value of $s$ is $k$.
(The function can be regarded as an instruction to differentiate
$x\cdot x \cdot x \cdots x$ first at factor number $s(1)$, then at
factor number $s(2)$, etc\., thus obtaining one of the
$m^i$ terms created by $i$ applications of the product rule.)
The left side of formula (b) is the same sum, taken over the
set of surjections.  The first term on the right side of (b) is
again this sum, now taken over the set of those
$s$ for which the restriction
to $\{1,\dots,i-1\}$ is surjective.
The second term is the same sum,
but taken over the set of those surjective
$s$ for which the restriction
to $\{1,\dots,i-1\}$ is not surjective.  From
these descriptions the equality is clear.
(This is essentially a proof of the basic recurrence formula for
Stirling numbers of the second kind;
cf\. formula (23), p\. 33 of [\stanley].)
\par
We prove (c) by induction on $i$.  When $i=0$ each term
of the sum except the first is divisible by $x$;
in this case the required polynomial in
$x, x', x'',\dots,x^{(n-j+1)}$ is zero.  For the inductive step,
apply $Q$ to both sides of the equation to obtain
$$
\align
Q^{i+1}(y^{(h)})
=&
\sum_{b=h}^{j-1}\frac{1}{(b-h)!} Q(y^{(b)})
	\lbrack Q^i(x^{b-h}) \rbrack_{x=0} \\
&+ \sum_{b=h}^{j-1}\frac{1}{(b-h)!}y^{(b)}
	Q\lbrack Q^i(x^{b-h}) \rbrack_{x=0} \\
&+ \text{ a polynomial in } x, x',\dots,x^{(n-j+1)}.
\endalign
$$
(Note that $Q(x^{(n-j+1)}) = 0$ by definition of
the differential operator
$Q$ in \thetag{\pandq}.)
In the first sum on the right, replace $b$ by $c-1$;
in all terms except the last replace
$Q(y^{(b)})$ by its chain rule equivalent $y^{(c)}x'$;
and absorb the last
term into the polynomial in $x, x',\dots,x^{(n-j+1)}$.
In each term of the second sum, use formula (b).
With these manipulations, we find that
$$
\align
Q^{i+1}(y^{(h)})
=&
\sum_{c=h+1}^{j-1}\frac{1}{(c-h-1)!} y^{(c)}x'
	\lbrack Q^i(x^{c-h-1}) \rbrack_{x=0} \\
&+ \sum_{b=h}^{j-1}\frac{1}{(b-h)!}y^{(b)}
	\left\{\lbrack Q^{i+1}(x^{b-h}) \rbrack_{x=0}
	-(b-h)x'\lbrack Q^i(x^{b-h-1}) \rbrack_{x=0} \right\} \\
&+ \text{ a polynomial in } x, x',\dots,x^{(n-j+1)} \\
=&
\sum_{b=h}^{j-1}\frac{1}{(b-h)!}y^{(b)}
	\lbrack Q^{i+1}(x^{b-h}) \rbrack_{x=0} \\
&+ \text{ a polynomial in } x, x',\dots,x^{(n-j+1)},
\endalign
$$
as required.
\par
To begin the proof of (d), observe that $P^h(x^b)$ vanishes if
$b<h$ and equals $\frac{b!}{(b-h)!}x^{b-h}$ otherwise.  Hence
$$
Q^i P^h(y) - \sum_{b = 0}^{j-1}\frac{1}{b!}y^{(b)}Q^i P^h(x^b)
=
Q^i(y^{(h)}) - \sum_{b=h}^{j-1}\frac{1}{(b-h)!}y^{(b)}Q^i(x^{b-h}).
\tag\base
$$
By statement (c), the right side of \thetag{\base} is,
modulo the principal ideal generated by $x$, a polynomial in
$x',x'',\dots,x^{(n-j+1)}$.  Since each term of the summation
involves some $y^{(b)}$, the polynomial in question is obtained by
expanding $Q^i(y^{(h)})$ and ignoring all terms
divisible by $x$ or any $y^{(b)}$.
Clearly all the coefficients of this polynomial are nonnegative.
\par
Let us define a  $\bold Z \oplus \bold Z$  grading
by declaring the degree of  $x^{(t)}$  to be  $(2 - t, 1)$
and the degree of  $y^{(t)}$  to be  $(2(j - t) - 1, j - 1 - t)$.
Then  $Q$  is a homomorphism of degree  $(-1, 0)$.
(Note in particular that  $y^{(j-1)}$  has degree $(1, 0)$
and that  $Q(y^{(j-1)}) = 1$.)
Each term in the expansion of $Q^i(y^{(h)})$ has degree
$$
( 2(j-1-h)+1-i,  j-1-h ).
$$
If  $i < 2(j - 1 - h) + 1$  then the first component is positive.
Hence each monomial in the expansion of $Q^i(y^{(h)})$
involves either  $x$  or  $x'$  or some $y^{(b)}$.
Hence the right side of \thetag{\base} is contained in the
ideal generated by $x$ and $x'$.
Similarly, if  $i = 2(j - 1 - h) + 1$  then the first component is zero.
Hence each monomial in the expansion of $Q^i(y^{(h)})$
either involves one of the same variables
or is a power of  $x''$;
the second component of the degree tells us that
the relevant power is  $(x'')^{j-1-h}$.
Clearly the coefficient of this term is not zero.
Hence the right side of \thetag{\base} is, modulo the
ideal generated by $x$ and $x'$,
a positive multiple of $(x'')^{j-1-h}$.
\par
The binomial formula for differential operators tells us that
the left side of formula (e) equals
$$
\sum_{a=0}^{i} \sum_{b=0}^{h} \binom{i}{a} \binom {h}{b}
Q^aP^b(\phi) \cdot Q^{i-a}P^{h-b}(x).
$$
Note that  $P^{h-b}(x)$  vanishes unless  $b$  equals  $h$  or  $h-1$;
that  $Q^{i-a}P(x)$ vanishes unless  $a=i$;
and that  $Q^{i-a}(x)= x^{(i-a)}$.  Hence
$$
Q^i P^h(x \cdot \phi)
= h Q^i P^{h-1}(\phi)
+ \sum_{a = 0}^{i}  \binom{i}{a} Q^a P^h(\phi) \cdot x^{(i-a)}.
$$
Modulo the ideal generated by $x$ and $x'$
the last two terms of the sum vanish.
\par
We prove (f) by induction on $k$.  The base case is statement (d).
For the inductive step, apply statement (e) to each
term on the left side of (f).
We find that, modulo the ideal generated by $x$ and $x'$,
$$
\align
Q^i P^h(x^ky) -& \sum_{b = 0}^{j-1} \frac{1}{b!}  y^{(b)} Q^i P^h(x^{k + b}) \\
& \equiv  h \lbrack Q^iP^{h-1}(x^{k-1}y)
	- \sum_{b = 0}^{j-1} \frac{1}{b!}  y^{(b)} Q^iP^{h-1}(x^{k+b-1}) \rbrack
 \\
&\quad + \sum_{a = 0}^{i-2}\binom{i}{a}\lbrack  Q^aP^h(x^{k-1}y)
	- \sum_{b = 0}^{j-1} \frac{1}{b!}  y^{(b)} Q^aP^h(x^{k+b-1}) \rbrack
	\cdot x^{(i-a)}.
\endalign
$$
Apply the inductive hypothesis to both bracketed expressions.
If  $i < 2(k + j - 1 - h) + 1$  then they both vanish.
If  $i = 2(k + j - 1 - h) + 1$  then the first term yields
a positive multiple of $(x'')^{k+j-1-h}$,
as does the last term ($a = i - 2$) in the sum.
Hence  $Q^i P^h(x^ky)$  is likewise a positive multiple of
$(x'')^{k+j-1-h}$.
\par
Statement (g) is the special case $h=j-1$ of statement (f).
\qed
\enddemo
\par
We now return to the proof of Lemma~A.
Recall that the first row of the matrix $M$
consists of all monomials of degree at most
$d$ in $x$ and $y$, and that a typical column
is shown in \thetag{\col}.
We wish to show that this matrix has
rank $n+1$ except possibly above the locus of
tangency to the divisor at infinity,
i\.e\., whenever $x' \neq 0$ or $x'' \neq 0$.
There are compatible actions of the projective
general linear group $PGL(2)$
on the universal family $\Cal C$
and on $F(1)$, the incidence correspondence of points and
lines in the plane.
Hence we may assume that we are studying a point
 of $\Cal C(n)$ over the point $x=y=y'=0$ of $F(1)$.
We may also assume that $x'=0$,
i\.e\., that the point lies over the
divisor at infinity, since we have already
dealt with other points when we
examined the primary charts.
\par
Consider the square submatrix of $M$ obtained by
extracting the columns headed by these
$n+1$ monomials:
$$
\gather
1, x, x^2, \dots, x^{j - 1}, \\
y, x^j, xy, x^{j + 1}, x^2y, x^{j + 2}, x^3y, \dots.
\endgather
$$
It  is of the form
$$
\bmatrix
A & C \\
B & D \\
\endbmatrix ,
$$
where the upper left  $j \times j$  submatrix  $A$
is upper triangular,
with diagonal entries $0!$, $1!$, $2!$, $3!$, $4!$, etc.
The lower left  $(n - j + 1) \times j$  submatrix  $B$  is a zero matrix.
The upper right  $j \times (n - j + 1)$  submatrix  $C$  is irrelevant.
We look at the lower right $(n - j + 1) \times (n - j + 1)$
matrix $D$
modulo the ideal generated by $x$ and $x'$.
Statement (a) of Lemma~B then says that
columns 2, 4, 6, etc\.---corresponding to the monomials
$x^j, x^{j + 1}, x^{j + 2},$ etc\.---are zero
above the diagonal, and have positive multiples
of powers of $x''$ on the diagonal.
Statement (g) of the same lemma
says that if one modifies
columns 1, 3, 5, etc\.---corresponding to the monomials
$y, xy, x^2y,$ etc\.---by
adding suitable linear combinations of previous columns,
then one again obtains, modulo the ideal,
columns which
are zero
above the diagonal, and have positive multiples
of powers of $x''$ on the diagonal.
Hence, modulo the ideal,
the submatrix  $D$  is
nonsingular whenever $x'' \neq 0$.
Hence the full square matrix is likewise nonzero,
and $M$ is of maximal rank $n+1$.
\qed
\enddemo
\par
We now consider
$\cn$, the fiber
product over $\pnd$ of the lifts of $\Cal C$
to $F(n_1)$, $F(n_2)$, $\dots$, $F(n_p)$.
We assume that each $n_i$ is positive.
Let $\fnopen$
denote the open subvariety of
$\fn=F(n_1)\times F(n_2)\times\dots\times F(n_p)$
obtained by removing the following points:
\roster
\item"$\bullet$" each point $(P_1,P_2,\dots,P_p)$
for which some $P_i$ lies on the intersection of
two divisors at infinity;
\item"$\bullet$" each point
for which some $P_i$ lies on a locus of tangency
to a divisor at infinity;
\item"$\bullet$" each point for which
some $P_i$ and some $P_j$
($i \neq j$)
lie over the same point of $\pp$
(i\.e\., each point lying over a large diagonal of
$\pp \times \pp \times\dots\times \pp$).
\endroster
Let
$\cnopen$
denote the inverse image of $\fnopen$ in
$\cn$.
\par
\proclaim{Lemma C}
Suppose that
$d \geq p-1 + \sum_{i=1}^p n_i$.
Then the morphism
$\cnopen \to \fnopen$
is smooth,  with relative dimension
$N(d) - \left(p + \sum_{i=1}^p n_i\right)$.
\endproclaim
Let $\fnplus$
denote the open subvariety of
$\fn$
obtained by removing the following points:
\roster
\item"$\bullet$" each point $(P_1,P_2,\dots,P_p)$
for which some $P_i$ lies on the intersection of
two divisors at infinity;
\item"$\bullet$" each point
for which some $P_i$ lies on a locus of tangency
to a divisor at infinity;
\item"$\bullet$" each point for which
some $P_i$ and some $P_j$
($i \neq j$)
lie over the same point of $F(1)$
(i\.e\., each point lying over a large diagonal of
$F(1)\times F(1)\times\dots\times F(1)$);
\item"$\bullet$" each point for which
some $P_i$, $P_j$, and $P_k$
($i$, $j$, $k$ distinct)
all lie over the same point of $\pp$;
\item"$\bullet$" each point for which
some $P_i$ and some $P_j$
($i \neq j$)
lie over the same point of $\pp$,
and one of them lies on a divisor
at infinity.
\endroster
Note that $\fnplus$ is larger than $\fnopen$,
since a point $(P_1,P_2,\dots,P_p)$ of
$\fnplus$ may project to a point
$(Q_1,Q_2,\dots,Q_p)$ of
$\pp \times \pp \times\dots\times \pp$
for which $Q_1$, $Q_2$, etc\., are not all distinct.
Let us stratify $\fnplus$ by
declaring a point $(P_1,P_2,\dots,P_p)$
to be in $\fnq$ if it projects to
a point $(Q_1,Q_2,\dots,Q_p)$ of
$\pp \times \pp \times\dots\times \pp$
for which the number of distinct points among
$Q_1$, $Q_2$, etc\., is exactly $q$.
The possible values for $q$ are
$\lceil \frac{p}{2} \rceil$, $\lceil \frac{p}{2} \rceil + 1$, \dots, $p$;
the stratum $\fn^p = \fnopen$
is open and dense.
Let
$\cnq$
denote the inverse image of $\fnq$ in
$\cn$.
\par
\proclaim{Lemma D}
Suppose that
$d \geq p-1 + \sum_{i=1}^p n_i$.
Then for each $q$ (from $\lceil \frac{p}{2} \rceil$ to $p$)
the morphism
$\cnq \to \fnq$
is smooth,  with relative dimension
$N(d) - \left(q + \sum_{i=1}^p n_i\right)$.
\endproclaim
\demo{Proof of Lemmas C and D}
Lemma C  is a special case of Lemma D.
\par
Over each point $(P_1,P_2,\dots,P_p)$ of $\fnq$
the fiber of $\cnq$ is a subspace of $\pnd$
of codimension at most $q + \sum_{i=1}^p n_i$,
as can be seen by a naive count of the number
of imposed conditions.
Hence to prove Lemma D it suffices to show that this fiber
is a linear subspace of
$\pnd$ of codimension at least $q + \sum_{i=1}^p n_i$.
This fiber is the intersection of
the fiber of $\Cal C(n_1)$ over $P_1$,
the fiber of $\Cal C(n_2)$ over $P_2$,
etc.
Lemma A therefore provides
an explicit set of generators
for its defining ideal
in $\pnd$.
These generators include, for each $P_i$ not
on a divisor at infinity, the functions
$$
f_i(x_i,y_i),\quad  f'_i(x_i,y_i,y'_i), \quad \dots \quad , \quad
f^{(n_i)}_i(x_i,y_i,y'_i,\dots,y^{(n_i)}_i)
\tag\generators
$$
(where $x_i$, $y_i$, etc\., are coordinates on
a primary chart of $F(n_i)$ containing $P_i$,
and $f_i$ generates the ideal of $\Cal C$)
together with, for each $P_i$
on a single divisor at infinity $I_{j_i}$, the functions
$$
\gather
f_i(x_i,y_i), \quad
Pf_i(x_i,y_i,y'_i), \quad \dots, \\
P^{j_i-1}f_i(x_i,y_i,y'_i, \dots, y^{(j_i-1)}_i), \quad
QP^{j_i-1}f_i(x_i,y_i,y'_i,\dots,y^{(j_i-1)}_i,x'_i), \\
\dots, \quad
Q^{n_i-j_i+1}P^{j_i-1}f_i(x_i,y_i,y'_i,\dots,y^{(j_i-1)}_i,x'_i,
	\dots,x^{(n_i-j_i+1)}_i)
\endgather
$$
(where $x_i$, $y_i$, etc\., are coordinates on
a secondary chart of $F(n_i)$ containing $P_i$).
Each of these generators is a linear function in
the $a_{uv}$'s; thus
the fiber of
$\cn$
over $(P_1,P_2,\dots,P_p)$
is defined by a system of $p + \sum n_i$ linear equations.
Let $M_{\bold n}$ denote the matrix of this system.
\par
Note that it suffices to prove the lemma
in case $d = p-1 + \sum_{i=1}^p n_i$, since increasing $d$ will only
enlarge $M_{\bold n}$ by adding more columns.
\par
We will use induction on $q$; the case $q=0$ is vacuous.
The inductive step has two cases.
Suppose first that the image of $P_1$ in $\pp$
is distinct from the images of $P_2,P_3,\dots,P_p$.
Let $\bold n - n_1$ denote $(n_2,n_3,\dots,n_p)$.
On $\pp$ choose an affine chart, with coordinates $x$ and $y$,
so that the chart contains the images of $P_2,P_3,\dots,P_p$
and so that the image of $P_1$ is the point at infinity on the
$y$-axis.
Around each $P_i$ ($i \neq 1$) choose, as appropriate,
a primary or secondary chart on $F(n_i)$
whose first two coordinates are $x_i=x$ and $y_i=y$.
(Choose a secondary chart only if the point lies on a divisor at infinity.)
Let $x_1=\frac{x}{y}$ and $y_1=\frac{1}{y}$;
then the image of $P_1$ is at the origin of the $\bold A^2$
with coordinates $x_1$ and $y_1$.
Note that over this affine chart the defining equation for
the universal family is
$$
\sum_{u+v \le d} a_{uv}{x_1}^u{y_1}^{d-(u+v)}=0.
\tag \aufeqn
$$
Around $P_1$ choose, as appropriate,
a primary or secondary chart on $F(n_1)$
whose first two coordinates are $x_1$ and $y_1$.
\par
We may assume that the matrix $M_{\bold n}$
has been arranged so that the initial columns correspond
to those $a_{uv}$'s for which $u+v \leq  d-1-n_1$.
Then $M_{\bold n}$ takes the form
$$
\bmatrix
M_{\bold n -n_1} & C \\
B & D \\
\endbmatrix .
$$
By the inductive hypothesis the rank of the upper left
submatrix is at least $q-1+\sum_{i=2}^{p} n_i$.
The $(1 + n_1)\times (1+N(d-1-n_1))$
submatrix $B$ records the partial derivatives
of the functions defining the fiber of
$\Cal C(n_1)$ over $P_1$
with respect to the $a_{uv}$'s just mentioned.
In \thetag{\aufeqn} the coefficient of each such $a_{uv}$
is divisible by ${y_1}^{n_1 +1}$.
If we implicitly differentiate this monomial
$n_1$ or fewer times, as we do when applying
the differential operator $P$ or $Q$,
the result is still divisible by $y_1$.
Hence all the entries of $B$ vanish
at $P_1$.
Submatrix $C$ is irrelevant.
Matrix $D$ contains the submatrix $M_{n_1}$;
by Lemma A its rank is $1 + n_1$.
Hence the rank of $M_{\bold n}$ is at least
$q + \sum_{i=1}^p n_i$.
\par
Next suppose that the image of $P_1$ in $\pp$
coincides with the image of some $P_i$ ($i \neq 1$).
Without loss of generality we may assume that
$P_1$ and $P_2$ project to the same point of $\pp$,
but that this point is distinct from
the images of $P_3,P_4,\dots,P_p$.
We may also assume that $n_1 \leq n_2$.
Furthermore the images of $P_1$ and $P_2$ in $F(1)$
are distinct points, and neither $P_1$ nor $P_2$
is on a divisor at infinity.
As in the previous case,
we choose on $\pp$ an affine chart, with coordinates $x$ and $y$,
so that the chart contains the images of $P_3,P_4,\dots,P_p$
and so that the common image of $P_1$ and $P_2$
is the point at infinity on the $y$-axis.
Around each $P_i$ ($i > 2$) we choose, as appropriate,
a primary or secondary chart on $F(n_i)$
whose first two coordinates are $x_i=x$ and $y_i=y$.
Let $x_1=x_2=\frac{x}{y}$ and $y_1=y_2=\frac{1}{y}$;
then the common image of $P_1$ and $P_2$
is at the origin of the $\bold A^2$
with coordinates $x_1$ and $y_1$.
Without loss of generality we may assume that
the image of $P_1$ in $F(1)$ represents the
tangent direction of the $y_1$-axis.
Hence we may use the primary chart with coordinates
$y_1$, $x_1$, $x'_1$, etc\.
on $F(n_1)$ and
the primary chart with
coordinates $x_2$, $y_2$, $y'_2$, etc\.
on $F(n_2)$.
\par
 We may assume that the matrix $M_{\bold n}$
has been arranged so that the initial columns correspond
to those $a_{uv}$'s for which $u+v \leq  d-2-(n_1+n_2)$.
Then $M_{\bold n}$ takes the form
$$
\bmatrix
M_{\bold n -n_1 -n_2} & C \\
B & D \\
\endbmatrix .
$$
By the inductive hypothesis the rank of the upper left
submatrix is at least $q-1+\sum_{i=3}^p n_i$.
The $(2 + n_1 + n_2)\times (1+N(d-2 - n_1 - n_2))$
submatrix $B$ records the partial derivatives
of the functions defining the fiber of
$\Cal C(n_1)$ over $P_1$
and the fiber of
$\Cal C(n_2)$ over $P_2$,
with respect to the $a_{uv}$'s just mentioned.
In \thetag{\aufeqn} the coefficient of each such $a_{uv}$
is divisible by ${y_1}^{n_1+n_2 +2}$.
If we implicitly differentiate this monomial
at most $n_2 = \max\{n_1,n_2\}$ times,
the result is still divisible by $y_1$.
Hence all the entries of $B$ vanish
at $(P_1,P_2)$.
Submatrix $C$ is again irrelevant.
Matrix $D$ contains the following submatrix:
$$
\bmatrix
1
   & x_1 & x_1^2 & \cdots & x_1^{n_2}
   & y_1 & y_1^2 & \cdots & y_1^{n_1}\\
0
   & x'_1 & 2x_1 x'_1  & \cdots & n_2x_1^{n_2-1}x'_1
   & 1 & 2y_1 & \cdots & n_1 y_1^{n_1-1}\\
\vdots
   & \vdots & \vdots & \ddots & \vdots
   & \vdots & \vdots & \ddots & \vdots \\
0
   & x_1^{(n_1)} &  & \cdots &
   & 0 & 0 & \cdots & n_1! \\
1
   & x_2 & x_2^2 & \cdots & x_2^{n_2}
   & y_2 & y_2^2 & \cdots & y_2^{n_1}\\
0
   & 1 & 2x_2  & \cdots & n_2x_2^{n_2-1}
   & y'_2 & 2y_2y'_2 & \cdots & n_1 y_2^{n_1-1}y'_2\\
\vdots
   & \vdots & \vdots & \ddots & \vdots
   & \vdots & \vdots & \ddots & \vdots \\
0
   & 0 & 0 & \cdots & n_2!
   & y_2^{(n_2)} &  & \cdots &
\endbmatrix .
$$
We can reorder the columns as follows:
$$
\bmatrix
1  & y_1   & x_1   & y_1^2   & x_1^2   & \cdots   & y_1^{n_1}  & x_1^{n_1}
& x_1^{n_1+1}  & \cdots  & x_1^{n_2}  \\
0  & 1   & x'_1   & 2y_1   & 2x_1 x'_1    & \cdots   & n_1 y_1^{n_1-1}
& \ddots  &   &   & \vdots  \\
0  & 0  & x''_1  & 2  & \ddots  &   &   &   &
&   & \\
\vdots  & \vdots   & \vdots   & \vdots  &    &    &
&   &   &   &   \\
0  & 0   & x_1^{(n_1)}  & 0  &  \cdots  &    &    &
&   &   & \\
1  & y_2   & x_2   & y_2^2   & x_2^2   & \cdots
& y_2^{n_1}  & x_2^{n_1}  & x_2^{n_1+1}
& \cdots  & x_2^{n_2}  \\
0  & y'_2   & 1   & 2y_2y'_2   & 2x_2    & \ddots
&  &   &   &   & \vdots  \\
0  & y''_2  & 0  & \ddots  &   &   &   &
&   &   &   \\
\vdots  & \vdots   & \vdots   &   &    &    &
&   &   &   &   \\
0  & y_2^{(n_2)}  & 0  & &   &   &
&   &   &   & \\
\endbmatrix .
$$
We can then rearrange the rows:
$$
\bmatrix
1  & y_1   & x_1   & y_1^2   & x_1^2   & \cdots   & y_1^{n_1}  & x_1^{n_1}
& x_1^{n_1+1}  & \cdots  & x_1^{n_2}  \\
1  & y_2   & x_2   & y_2^2   & x_2^2   & \cdots
& y_2^{n_1}  & x_2^{n_1}  & x_2^{n_1+1}
& \cdots  & x_2^{n_2}  \\
0  & 1   & x'_1   & 2y_1   & 2x_1 x'_1    & \cdots   & n_1 y_1^{n_1-1}
& \ddots  &   &   & \vdots  \\
0  & y'_2   & 1   & 2y_2y'_2   & 2x_2    & \ddots
&  &   &   &   &   \\
0  & 0  & x''_1  & 2  & \ddots  &   &   &   &
&   & \\
0  & y''_2  & 0  & \ddots  &   &   &   &
&   &   &   \\
0  & 0  & x^{(3)}_1  &   &   &   &   &
&   &   &   \\
\vdots  & \vdots   & \vdots   &   &    &    &
&   &   &   &   \\
0  & 0   & x_1^{(n_1)}  &   &    &    &    &
&   &   &   \\
0  & y_2^{(n_1)}  & 0  &   &   &   &   &
&   &   &   \\
0  & y_2^{(n_1+1)}  & 0  &   &   &   &
&   &   &   &   \\
\vdots  & \vdots  & \vdots  &   &   &   &
&   &   &   &   \\
0  & y_2^{(n_2)}  & 0  & \cdots  &   &   &
&   &   &   &
\endbmatrix .
$$
Evaluating at $x_1 = y_1 = x'_1 = x_2 = y_2  = 0$, we obtain
$$
\bmatrix
1 & 0 & 0 & 0 & 0 & \cdots & 0 & 0 & 0 & 0 & \cdots & 0\\
1 & 0 & 0 & 0 & 0 & \cdots & 0 & 0 & 0 & 0 & \cdots & 0\\
& 1 & 0 & 0 & 0 & \cdots & 0 & 0 & 0 & 0 & \cdots & 0\\
&  & 1 & 0 & 0 & \cdots & 0 & 0 & 0 & 0 & \cdots & 0\\
&  &  & 2! & 0 & \cdots & 0 & 0 & 0 & 0 & \cdots & 0\\
&  &  &  & 2! &  \cdots & 0 & 0 & 0 & 0 & \cdots & 0\\
& & & & & \ddots & \vdots & \vdots & \vdots & \vdots & \ddots & \vdots\\
& & & & & & n_1! & 0 & 0 & 0 & \cdots & 0\\
& & & & & & & n_1!  & 0 & 0 & \cdots & 0\\
& & & & & & & & (n_1+1)!  & 0 & \cdots & 0\\
& & & & & & & & & (n_1+2)!  & \cdots & 0\\
& & & & & & & & & & \ddots & \vdots\\
& & & & & & & & & & & n_2!\\
\endbmatrix ,
$$
which clearly has rank
$n_1+n_2+1$.
Hence at $(P_1,P_2)$ the rank of $D$
is at least $n_1+n_2+1$, and
the rank of $M_{\bold n}$ is at least
$q + \sum_{i=1}^p n_i$.
\qed
\enddemo
\par
By definition $\Cal C(n)$ is the union of
$\lambda(\Cal C)$ (where $\lambda$
is the lifting map $\Cal C \dashrightarrow  F(n)$) and a closed subvariety
$\Cal S(n)$, each point of which lies
over a singular or nonreduced point of some member of $\Cal C$.
Since the codimension of
$\Cal C(n)$ in $\bold A^{n+2} \times \pnd$
is $n+1$, the codimension of $\Cal S(n)$ is at least $n+2$.
\proclaim{Lemma E}
Suppose that $d \geq n$.
Suppose that $P$ is a point of $F(n)$ which is not on any
divisor at infinity.
Then the fiber of $\Cal S(n)$ over $P$
has codimension at least $n+2$ in $\pnd$.
\endproclaim
\demo{Proof}
It suffices to check over the primary chart of $F(n)$ with
coordinates $x, y, y', y''$, $y^{(3)},\dots, y^{(n)}$.
Over this chart the defining ideal of $\Cal S(n)$  contains
the functions $f, f', f''$, $f^{(3)},\dots, f^{(n)}$
defining $\Cal C(n)$
(obtained from the function $f$ of (\ufeqn)
by repeated implicit differentiation)
and the partial
derivative $\frac{\partial f}{\partial y}$.
Let $M$ denote the $(n+2)\times({N(d)}+1)$
matrix of partial derivatives of
these functions
with respect to each of the
$a_{uv}$.
(Once again we may interpret $M$ as the matrix of a linear system.)
The square submatrix consisting of the partials with respect to
$a_{00}, a_{10}, a_{20},\dots,a_{n0}$, and $a_{01}$
is nonsingular at each point:
$$
\bmatrix
1 & x & x^2 & \cdots & x^n & y \\
0 & 1 & 2x  & \cdots & nx^{n-1} & y' \\
0 & 0 & 2    & \cdots & n(n-1)x^{n-2} & y'' \\
\vdots & \vdots & \vdots & \ddots & \vdots & \vdots \\
0 & 0 & 0 & \cdots & n! & y^{(n)} \\
0 & 0 & 0 & \cdots & 0 & 1
\endbmatrix .
\tag\matparttwo
$$
Hence at each point the functions
$f, f', f'', f^{(3)},\dots,f^{(n)}$, and $\frac{\partial f}{\partial y}$
define a variety of codimension $n+2$
which is smooth over the primary chart of $F(n)$.
\qed
\enddemo
\par
Let $\sn$
denote the fiber product
$$
\Cal S(n_1) \fpopn \Cal C(n_2)\fpopn \cdots \fpopn \Cal C(n_p).
$$
Let $\fn^{\text{finite}}$ denote the open subvariety of $\fn$
obtained by removing, from each factor, all divisors at infinity;
let $\sn^{\text{finite}}$ denote its inverse image in $\sn$.
Let
$$
\snopen = \cnopen \cap \sn^{\text{finite}}.
$$
\par
\proclaim{Lemma F}
Suppose that
$d \geq p-1 + \sum_{i=1}^p n_i$.
Then each fiber of the morphism
$\snopen \to \fnopen \cap \fn^{\text{finite}}$
has codimension at least $1 + p + \sum_{i=1}^p n_i$.
\endproclaim
\demo{Proof}
The proof is similar to that of Lemma D.
In this instance, we need only to work in primary charts.
The defining ideal for the fiber of $\sn$
over a point $(P_1,P_2,\dots,P_p)$
of $\fn^{\text{finite}}$
includes the functions of \thetag{\generators}
and the partial derivative
$\frac{\partial f_1}{\partial y_1}(x_1,y_1)$.
Each of these functions is linear in
the $a_{uv}$'s; let
$M$ denote the matrix of this linear system.
As in the proof of Lemma D,
we note that it suffices to prove the lemma
in case $d = p-1 + \sum_{i=1}^p n_i$.
\par
If $(P_1,P_2,\dots,P_p)$ is in $\snopen$,
then the images in $\pp$ of these $p$ points are all
distinct.
By Lemma D, the fiber of
$\Cal C_{\pnd}(n_2,n_3,\dots,n_p)$
over $(P_2,P_3,\dots,P_p)$
has codimension $p - 1 + \sum_{i=2}^p n_i$.
By Lemma E, the fiber of
$\Cal S(n_1)$ over $P_1$
has codimension at least $n_1+2$.
Since the image of $P_1$ is distinct from the
images of $P_2$ through $P_p$,
we may argue as in the first case
of the inductive step in the proof of Lemma D.
By an appropriate choice of affine charts,
in particular, by choosing a chart on $\pp$
containing the images of $P_2,P_3,\dots,P_p$ and
so that the image of $P_1$
is the point at infinity on the $y$-axis,
we see that $M$
takes the form
$$
\bmatrix
M_{\bold n -n_1} & C \\
B & D \\
\endbmatrix ,
$$
where $M_{\bold n -n_1}$
has rank $p - 1 + \sum_{i=2}^p n_i$,
$B$ is a zero matrix, and
$D$ has rank at least $n_1 + 2$.
Hence $M$ has rank at least $1 + p + \sum_{i=1}^p n_i$.
\qed
\enddemo
\par
\demo{Proof of Theorem 3}
Let
$\bold n=(n_1,n_2,\dots,n_p)=(o_1-1,o_2-1,\dots,o_p-1)$.
By Lemma D,
each morphism $\cnq \to \fnq$ is smooth.
\par
Let us assume for the moment that $S$ is nonsingular.
Let $\tau:S \to \pnd$ be the morphism determined by the
family $\Cal X$.  Then the morphism
$$
\gathered
S \times PGL({N(d)}) \to \pnd \\
(s,\gamma) \mapsto \gamma \cdot \tau(s)
\endgathered
\tag\gp
$$
is smooth.
(This assertion is justified in the course of the proof of
Theorem~2 in [\ktrans].)
Hence the projection
$$
\cnq \fpopn
(S \times PGL({N(d)})) \to \cnq,
$$
which is obtained from \thetag{\gp} by base extension, is smooth.
Composing it with the smooth morphism
$\cnq \to \fnq$,
we obtain
$$
\cnq
\fpopn (S \times PGL({N(d)})) \to \fnq,
$$
which fits into the following diagram.
(The vertical morphism on the left is the composite of two projections;
the one on the right is inclusion.)
$$
\CD
@.  C_1(n_1)\times  \dots \times C_p(n_p) \cap  \fnq
	\negthickspace \negthickspace
	\negthickspace \negthickspace
	\negthickspace \negthickspace
	\negthickspace \negthickspace
	\negthickspace \negthickspace
	\negthickspace \negthickspace
	\negthickspace \negthickspace
	\negthickspace \negthickspace
	\negthickspace \negthickspace
	\negthickspace \negthickspace \\
@. @VVV @. \\
\cnq
	\fpopn (S \times PGL({N(d)})) @>>>
	\fnq \thickspace \\
@VVV @. \\
PGL({N(d)}) @.
\endCD
$$
\par
The fiber of
$\cnq
\fpopn (S \times PGL({N(d)}))$
over a point
$\gamma \in PGL({N(d)})$
is the fiber product
$\cnq \fpopn S$, where $S$ is
regarded as a $\pnd$-variety
via the morphism $\gamma \circ \tau$, i\.e\.,
via $\tau$ followed by translation by $\gamma$.
Kleiman's transversality lemma
(Lemma~1 of [\ktrans]) tells us that for generic $\gamma$
the map
$\cnq \fpopn S \to \fnq$
is transverse to
$C_1(n_1)\times  \dots \times C_p(n_p) \cap  \fnq$.
Our assumption that the family $\Cal X$ is generic
means that we may assume $\gamma$ is the identity.
\par
We now count dimensions.  By Lemma D,
$$
\dim \left( \cnq \right)
	= \dim \fnq + N(d) - \left(q + \sum_{i=1}^p n_i\right).
$$
Since the morphism $\tau$ is generic,
$$
\dim \left( \cnq \fpopn S \right)
	= \dim \fnq  - \left(q + \sum_{i=1}^p n_i\right) + s = \dim \fnq -q.
$$
To specify a point of $C_1(n_1)\times  \dots \times C_p(n_p) \cap  \fnq$,
one can first specify $p-q$ points
lying on the intersection of two of the $p$ curves
$C_1$, \dots, $C_p$,
then specify $2q-p$ additional points, each of which lies on one
of the curves; there may be further choices involved in specifying
higher-order data,
but these choices cannot contribute to the dimension count:
$$
\dim \left( C_1(n_1)\times  \dots \times C_p(n_p) \cap  \fnq \right)
	= 2q - p .
$$
If $q < p$, the dimension count shows that the image of
$\cnq \fpopn S$ is disjoint from
$C_1(n_1)\times  \dots \times C_p(n_p)$.
Hence the morphism
$\cnplus \fpopn S \to \fnplus$
is transverse to
$C_1(n_1)\times  \dots \times C_p(n_p) \cap \fnplus$,
and all the intersections lie in the open dense stratum $\fnopen$.
Our assumptions on the $p$ curves guarantee that
$C_1(n_1)\times  \dots \times C_p(n_p)$
is contained in $\fnplus$.  Hence
enlarging $\cnplus$ to $\cn$
creates no further intersections:
the map from
$\cn \fpopn S$
to  $\fn$
is likewise transverse to the $p$-fold product.
\par
Note that $\Cal X_S(\bold n)$,
the fiber product
over $S$ of the various lifts of $\Cal X$,
is a subvariety of
$\cn \fpopn S$,
and that it has the same dimension.
Hence the morphism
$\sigma: \Cal X_S(\bold n)  \to \fn$
(the restriction of the projection of $\fn \times S$
onto its first factor)
is also transverse to
$C_1(n_1)\times C_2(n_2)\times\dots\times C_p(n_p)$.
Each intersection between $\sigma(\Cal X_S(\bold n))$
and the $p$-fold product
is a $p$-tuple
$(x_1,x_2,\dots,x_p)$
in which $x_1$ is a contact or a false contact between $C_1$
and some member of $\Cal X$, and $x_2$ is a contact or a false contact
between $C_2$
and the same member of $\Cal X$, etc.
The number of such intersections is
$$
\int_{\fn}
	\lbrack
	C_1(n_1)\times C_2(n_2)\times \dots \times C_p(n_p)
	\rbrack
	\cdot
	\sigma_*\lbrack \Cal X_S(\bold n) \rbrack ,
$$
which is equal to the proto-contact number as defined by
\thetag{\pcn} in \S 3.
\par
We must show that each point of intersection is a $p$-tuple
$(x_1,x_2,\dots,x_p)$
of honest contacts rather than false ones.
Consider  $C_1(n_1)_{\text{sing}}$, the (zero@-dimensional) subvariety
of $C_1(n_1)$ lying over the singularities of $C_1$.
The dimension of
$C_1(n_1)_{\text{sing}}\times C_2(n_2)\times \dots \times C_p(n_p)$
is $p-1$.  Hence by Kleiman's transversality lemma the image
of $\cn \fpopn S$
in $\fn$ is disjoint from this $p$-fold product.
{\it A fortiori}, the image of $\Cal X_S(\bold n)$
is disjoint from this $p$-fold product.
Hence $x_1$ does not lie over a singular point of $C_1$.
Note this implies that $x_1$ does not lie over a divisor
at infinity on  $F(n_1)$.
Similarly one sees that $x_2$ does not lie over a singular point of $C_2$, and
hence that $x_2$ does not lie over a divisor
at infinity on  $F(n_2)$, etc.
Lemma F, together with the appropriate dimension count, shows that
the image of $\snopen \fpopn S$
is disjoint from $C_1(n_1)\times\dots\times~C_p(n_p)$.
As we have already seen, enlarging $\snopen$ to $\sn^{\text{finite}}$
creates no intersections.
Hence for every point of intersection $(x_1,x_2,\dots,x_p)$
between the image of $\Cal X_S(\bold n)$
and the $p$-fold product, $x_1$ lies over a nonsingular point
of the relevant member of $\Cal X$.
Similarly one sees that $x_2$ lies over a nonsingular point
of the relevant member of $\Cal X$, etc.
\par
If $S$ is singular, we can apply our argument above to
the singular locus $S_{\text{sing}}$.
The morphism
$$
S_{\text{sing}} \times PGL({N(d)}) \to \pnd
$$
is now flat rather than smooth.  (Again, this assertion is justified by
Kleiman in the proof of his Theorem~2.)
Kleiman's lemma now tells us, after the appropriate dimension count,
that the image of
$\Cal X_{S_{\text{sing}}}(\bold n)$
in $\fn$ is disjoint from
$C_1(n_1)\times\dots \times  C_p(n_p)$.
Hence, for a generic family $\Cal X$,
the members over  $S_{\text{sing}}$
have no simultaneous contacts with $C_1, C_2, \dots, C_p$.
Thus we may assume, as we did, that $S$ is nonsingular.
\qed
\enddemo
\par
\bigpagebreak
\heading 5. The higher-order characteristic numbers
of a family of plane curves
\endheading
\smallpagebreak
\par
At the beginning of \S 3
we defined the higher-order characteristic numbers
of a plane curve.  The degree and class are well known.
Each of the other numbers counts, perhaps with multiplicities,
the number of cusps of a specified order.  These numbers are
readily calculated, either implicitly or from a parametrization
of the curve; in fact one needs only a local or even formal
parametrization at each singular point.
(See the algorithm presented in the proof of
Proposition 3\.9 of [\ckho].)
\par
For a family of curves,
the numbers defined by \thetag{\indcnf} are
likewise called
{\it characteristic numbers}.
For an $s$-parameter family $\Cal X$, this formula
associates such a number to each monomial of weight $s$
in the indeterminates of \thetag{\indet}.
To denote this number we use the monomial itself,
with lowercase Greek letters replacing their
uppercase counterparts.
\par
Note that $\Lambda_0$ has weight 0.  If  $M$ is a monomial
of weight $s$ with associated characteristic number $m$, then
it is easy to show that
the characteristic number associated to
$\Lambda_0 M$ is $dm$, where $d$ is the degree of the
general member of $\Cal X$.
In the remainder of this section we consider
monomials not involving $\Lambda_0$.
\par
If $M$ is a monomial in $\Lambda_1$ and $\Pi_1$,
then we call its associated characteristic number {\it ordinary}.
It is well known that, under mild hypotheses on $\Cal X$,
the characteristic number
$(\lambda_1)^r (\pi_1)^{s-r}$ is the number of members
tangent to $r$ specified general lines and passing through $s-r$
specified general points.
These characteristic numbers, especially those of plane cubics,
have been the subject of numerous
investigations, both classical and contemporary.
(For example, see [\aone], [\atwo], [\athree],
[\ksone], [\kstwo], [\ksthree],
[\maillard], [\mx], [\schub], [\zeuthen].)
If $M$ also involves the indeterminate $\Gamma^2_2$,
then it should perhaps still be called ``ordinary'',
since its definition uses
only the notion of ordinary contact.  For example if $\Cal X$
is a general two-parameter family then
$$
\gamma^2_2
  = \int_{F(1)} \overline {\Gamma^2_2} \cap \sigma_* [\Cal X(1)]
  = \int_{F(1)} h^2 \check h \cap \sigma_* [\Cal X(1)]
$$
is the number of members
tangent to a specified line at a specified point.
This characteristic number does not, however, appear in
classical ordinary contact formulas.  It does appear in Schubert's
formula for triple contacts between a two-parameter family
and a single specified curve.  (See [\cktrip], [\rstwo], [\schubtwo],
and \S 6(a) of the present paper.)
\par
If $M$ involves other indeterminates, then
the associated number $m$ is called a
{\it higher-order characteristic number}.
For example, a two-parameter family has
six characteristic numbers.  The ordinary characteristic
numbers are
$(\lambda_1)^2$, $\lambda_1\pi_1$, $(\pi_1)^2$,
and $\gamma^2_2$.
The higher-order characteristic numbers are
$$
\align
\lambda_2 &= \int_{F(2)} \overline \Lambda_2 \cap \sigma_* [\Cal X(2)]
	= \int_{F(2)} \check h^2 z_2 \cap \sigma_* [\Cal X(2)];\\
\text{and }
\quad \pi_2 &= \int_{F(2)} \overline \Pi_2 \cap\sigma_* [\Cal X(2)]
	= \int_{F(2)} h^2 i_2 \cap \sigma_* [\Cal X(2)].
\endalign
$$
Theorem~3 tells us that, for a generic family of curves of
degree $\geq 2$,
$\lambda_2$ is the number of triple contacts
between a member of $\Cal X$ and a specified line.
(More generally, $\lambda_s$ is,
for a generic $s$-parameter family $\Cal X$
of curves of degree at least $s$,
the number of contacts of order $s+1$
between a member of $\Cal X$ and a specified line.)
One can show that $\pi_2$ is the number
of members of $\Cal X$ having a cusp at a specified point
(See [\cktrip], Theorem~A1.  Note that a ``dual" characteristic
number involving flexes rather than cusps
is considered in [\schubtwo] and [\rstwo].)
\par
The definition
$$
\pi_s  = \int_{F(s)} h^2 i_2 i_3 \cdots i_s \cap \sigma_* [\Cal X(s)]
$$
suggests that this characteristic number measures,
for a generic $s$-parameter family $\Cal X$
of curves of sufficiently large degree,
the number of member curves of $\Cal X$ whose lift,
at a specified point, meets the divisors at infinity
$I_2$, $I_3$, \dots, $I_s$.  (One might call such a point
a ``super-profound cusp''.)
To justify this interpretation would involve an appeal
to Kleiman's transversality lemma, to rule
out contributions to the intersection number created by the
way in which curvilinear data specializes at
a singular or nonreduced member of the family.
But our crucial Lemma A fails in precisely this
``super-profound'' situation,
and we cannot use Kleiman's lemma.
We suspect that the naive interpretation of  $\pi_s$ is incorrect,
and that one cannot even guarantee that the intersection
in question can be made proper.
\par
Similar difficulties plague the interpretation of other characteristic
numbers.   The number of characteristic numbers also grows
rather quickly with the number of parameters.
For example, a five-parameter family
possesses 70 characteristic numbers.
\par
\bigpagebreak
\heading 6. Variations, further remarks, and an example.
\endheading
\smallpagebreak
\par
\subheading{(a) Contacts between a family and a single curve}
\smallpagebreak
\par
When $p = 1$, the special case of the contact formula says that
the number of contacts of order $o = n+1$ between
a specified curve $C$ (with no profound cusp)
and a generic $n$-parameter family $\Cal X$ is
$$\gather
\int_{F(n)} [C(n)]\cap \sigma_* [\Cal X(n)] \\
= d\lambda_n + \check d\pi_n + (3\check d + \kappa_2)\gamma_n^2
+ \dots
+ ((n+1)\check d + n\kappa_2 + \dots
+ 3\kappa_{n-1} + \kappa_n)\gamma_n^n
\endgather$$
where
$$\align
\lambda_n &= \int \check h^2 z_n \cap \sigma_*[\Cal X(n)],\\
\pi_n &= \int h^2 i_2 \cdots i_n \cap \sigma_*[\Cal X(n)],\\
\gamma_n^k &=
\int h^2 \check h z_{k-1}i_{k+1}i_{k+2} \cdots i_n
\cap \sigma_*[\Cal X(n)], \quad
k = 2,\dots,n.
\endalign$$
Note that this formula is essentially implicit in
the definition of the contact module $m_n(C)$.
\subheading{(b) Simultaneous contact with nonsingular curves}
\smallpagebreak
\par
Suppose that the curves $C_1,\dots, \allowmathbreak C_p$
are all nonsingular.  Then for each curve the characteristic numbers
$\kappa_1, \kappa_2, \dots$
vanish and (by the
Pl\"ucker formula) the class
is $\check d = d (d - 1)$.
In this case the contact formula tells us that
the number of simultaneous contacts of order
$(o_1,\dots,o_p) = (n_1+1,\dots,n_p+1)$
between $C_1,\dots,C_p$ and some member of a
$\sum n_j$@-parameter family
$\Cal X$
is obtained from
$$
\prod_{j=1}^p \medspace [d_j \Lambda_{n_j} +
d_j (d_j - 1)(\Pi_{n_j} + \sum_{k=2}^{n_j}(k+1) \Gamma_{n_j}^k)]
\tag \nonsing
$$
by the expansion and evaluation of monomials described
in Theorem~2.
\par
This contact formula
is valid under the assumptions of Theorem~3, that is,
when $\Cal X$ is a generic family of plane curves of
sufficiently high degree.
There is, however, a certain freedom in choosing hypotheses.
We could assume, for example, that each individual curve
$C_j$ is generic (hence nonsingular) of degree $d_j \ge o_j - 1$,
but make no assumption whatever about the family
(except, of course,
that the generic member is a reduced curve).
Then, if the base field
is either of characteristic zero or of
characteristic at least $\max\{o_1,\dots,o_p\}$,
formula \thetag{\nonsing}
counts (again, after expansion and evaluation of monomials)
the number
of simultaneous contacts.  Moreover, if each $d_j \ge o_j$,
then all contacts are of order exactly $(o_1,\dots,o_p)$
(that is, none of the
contacts with any of the $C_j$'s is of higher order).
The proof is essentially that of
Proposition 2\.5 of [\ckho], which treats the case $p=1$.
\par
\subheading{(c) B\'ezout's Theorem}
\smallpagebreak
\par
As we remarked in \S 4, Theorem~3
is valid even if we drop the requirement that each
specified order of contact $o_i$ be at least 2.
Suppose, for example, that $C$ and $C_1,C_2,\dots,C_p$
are reduced plane curves of degrees $d, d_1, d_2, \dots, d_p$
respectively.
Then our contact formula says that
the number of simultaneous intersections between
$C$ and $C_1,C_2,\dots,C_p$ is
$d^p d_1 d_2 \dotsm d_p$.
This is just a silly formulation of B\'ezout's Theorem.
Indeed, a simultaneous intersection is a $p$-tuple
$(x_1,x_2,\dots,x_p)$
in which the point $x_1$ lies in the intersection of $C$ and $C_1$,
the point $x_2$ lies in the intersection of $C$ and $C_2$, etc.
There are $d d_1$ possibilities for $x_1$, together with
$d d_2$ possibilities for $x_2$, etc.
The other cases omitted from our statement of Theorem~3
are of the same ilk.
\par
\subheading{(d) The formula of Fulton, Kleiman, and MacPherson}
\smallpagebreak
\par
The overlap of our contact formula and that of
Fulton-Kleiman-MacPherson
[\fkm] is the following formula:
\par
\proclaim\nofrills{}
The number of simultaneous ordinary contacts between
the members of a $p$-parameter family $\Cal X$ of plane curves
and $p$ specified plane curves $C_1$
(of degree $d_1$ and class $\check d_1$),
\dots, $C_p$ (of degree $d_p$ and class $\check d_p$)
is obtained from the product of modules
$$
\prod_{j=1}^p (d_j \Lambda + \check d_j \Pi)
$$
by expansion and evaluation of monomials.  To evaluate
$\Lambda^r\Pi^{p-r}$ means to replace it by the (ordinary)
characteristic number $\lambda^r\pi^{p-r}$,
the number of members of $\Cal X$
tangent to $r$ specified general lines and passing through $p-r$
specified general points.
\endproclaim
\par
Fulton {\it et al\.} assume that the individual curves
$C_1, C_2,\dots, C_p$ are in general position,
whereas our hypotheses involve genericity assumptions
about the family $\Cal X$ of curves, as well as different,
but mild, conditions on $C_1,\dots, C_p$.
(See (f) below for further comments along these lines.)
\par
\subheading{(e) The formula of de~Jonqui\`eres}
\smallpagebreak
\par
The classical formula of de~Jonqui\` eres [\dj] gives
the number of plane curves
of degree $d$ making contacts of
orders $o_1,o_2,\dots,o_p$ with a given curve
$C$ and which pass through $k$ points,
where $k = d(d+3)/2 + p - \sum_j o_j$.
It is tempting to use our set-up to obtain this formula
by calculating
$$
I:=\int_{\fn}
	\pi_1^*(h^2)
	\cdot
	\dots
	\cdot
	\pi_k^*(h^2)
	\cdot
	\pi_{k+1}^*\lbrack C(n_1)\rbrack
	\cdot
	\dots
	\cdot
	\pi_{k+p}^*\lbrack C(n_p)\rbrack
	\cap
	\sigma_*\lbrack \cn \rbrack .
$$
where now
$$
\gather
\fn := \undersetbrace \text{$k$ factors}
\to{\pp \times \dots  \times\pp}
 \times F(n_1) \times \dots \times F(n_p),\\
\cn := \undersetbrace \text{$k$ factors}
\to{\Cal C \fpopn \cdots \fpopn \Cal C} \fpopn
\Cal C(n_1)\fpopn
\cdots \fpopn \Cal C(n_p),
\endgather
$$
each $n_j=o_j -1$,
and $\Cal C$ denotes the universal family of degree $d$
plane curves.
However, Theorem~3 itself will never apply to
the universal family.
We have $\dim \cn = N(d) + p + k$ and
$\dim \fn = 2k + p + \sum_j o_j$ so that the map
$\cn \to \fn$ is smooth only if
$N(d) - k > \sum_j o_j$ while we must have $N(d) - k = \sum_j o_j - p$
for the proto-contact formula to be valid.
\par
Thus our formula represents an approach to problems of
higher-order contact distinct from that provided by the
formula of
de~Jonqui\` eres.  Note that de~Jonqui\` eres's formula
can lose its enumerative significance when,
for example, members of the family have nonreduced
components of sufficiently high multiplicity, whereas
the hypotheses of Theorem~3 guarantee validity of
our formula for suitable families of plane curves of fixed degree.
(See [\fkm] p\. 184ff for further
discussion of the enumerative significance
of de~Jonqui\` eres's formula.)
Moreover,  our results are
in keeping with the spirit of Hilbert's 15th problem
 which asks ``to establish rigorously and with an
exact determination of the limits
of their validity those geometrical numbers \dots" [\hilbert].
Finally, Fulton, Kleiman, and MacPherson [\fkm] note that
``[i]t was observed long ago that
de~Jonqui\` eres's formula yields via symbolic
multiplication a formula
in the case of several given curves."
The proto-contact formula
of Theorem~2 involves just such a general multiplication.
\par
Other contemporary treatments of de~Jonqui\` eres's formula
include [\mattuck] and [\vain].  In these cases,
the problem is formulated in terms
of systems of divisors of degree $d$ cut out on
the fixed curve $C$ by a
family $\Cal X$ of curves.  The formula so produced
counts the number of
such divisors on $C$ that have coalesced
in a prescribed manner.  In particular,
the support of such degenerate divisors could include points
that are singular either on $C$ or on the relevant
member of $\Cal X$.
In our language, such divisors
represent ``false contacts" and are not counted by our formula.
(See Example 3\.8, pp\. 503--504 of [\ckho] for
an illustration of this phenomenon.)
\par
\subheading{(f) An example: triple contact with two curves}
\smallpagebreak
\par
In this case we have two specified curves $C$ and $D$,
and seek the number of simultaneous contacts of order $(3,3)$
with members of a 4-parameter family $\Cal X$.
If $C$ and $D$ intersect transversely and have
no profound or flat cusps,
and if $\Cal X$ is
a generic family of curves of degree $d \ge 5$, then
Theorems~2 and~3 show that the desired number
$$
I=\int_{F(2) \times F(2)}
	\pi_1^*\lbrack C(2)\rbrack
	\cdot
	\pi_2^*\lbrack D(2)\rbrack
	\cap
	\sigma_*\lbrack \Cal X_S(2,2) \rbrack
$$
is obtained from the product of modules
$$
(d_C \Lambda_2 + \check d_C \Pi_2 +
	(3\check d_C + \kappa_{2C})\Gamma_2^2)
(d_D \Lambda_2 + \check d_D \Pi_2 +
	(3\check d_D + \kappa_{2D})\Gamma_2^2),
$$
where $d_C$, $\check d_C$, $\kappa_{2C}$
are the characteristic numbers of $C$, etc.
Explicitly, the number of simultaneous contacts is
$$
\gathered
d_C d_D(\lambda_2)^2
+ (d_C \check d_D +d_D \check d_C)\lambda_2 \pi_2 \\
\quad +
	(d_C(3\check d_D + \kappa_{2D}) +
		d_D(3\check d_C + \kappa_{2C}))
	\lambda_2 \gamma_2^2
+ \check d_C \check d_D(\pi_2)^2 \\
\quad +
		(\check d_C(3\check d_D + \kappa_{2D}) +
		\check d_D(3\check d_C + \kappa_{2C}))
	\pi_2 \gamma_2^2
+ (3\check d_C + \kappa_{2C})(3\check d_D + \kappa_{2D})
	(\gamma_2^2)^2 .
\endgathered
\tag\triptwo
$$
\par
Another approach can be taken to establishing
the validity of this
contact formula.  In [\cktrip] and [\rstwo], the
$PGL(2)$-orbits of $F(2)$ are identified.  There are three
of them:
the dense orbit $\Cal O(\dash)$, each point of which is
represented by the germ of a nonsingular curve without a flex,
the special orbit $\Cal O(0) = Z_2$, each point of which is
represented by the germ of a line,
and the divisor at infinity $\Cal O(\infty) = I_2$, each
point of which is
represented by the germ of an ordinary cusp.
Thus the action of $PGL(2) \times PGL(2)$
on $F(2) \times F(2)$
has nine orbits:
\blank
{\settabs 2 \columns
\+ \qquad \qquad \qquad {\bf Orbit} & {\bf Dimension}  \cr
\+ \qquad \qquad \qquad $\Cal O(\dash) \times \Cal O(\dash)$
	& \qquad 8 \cr
\+ \qquad \qquad \qquad $\Cal O(\dash) \times \Cal O(0), \quad
	\Cal O(0) \times \Cal O(\dash)$ &  \qquad 7 \cr
\+ \qquad \qquad \qquad $\Cal O(\dash) \times \Cal O(\infty), \quad
	\Cal O(\infty) \times \Cal O(\dash)$ & \qquad 7 \cr
\+ \qquad \qquad \qquad $\Cal O(0) \times \Cal O(0), \quad
	\Cal O(\infty) \times \Cal O(\infty)$ &  \qquad 6 \cr
\+ \qquad \qquad \qquad $\Cal O(0) \times \Cal O(\infty), \quad
	\Cal O(\infty) \times \Cal O(0)$ &  \qquad 6 \cr}
\blank
If  $C$, $D$, and $\Cal X$ are suitably transverse,
i\.e\., if $\Cal X_S(2,2)$ is mapped transversely to
$C(2) \times D(2)$ by $\sigma$ and the intersection
$(C(2) \times D(2))\cap \sigma(\Cal X_S(2,2)) $
is contained in the dense orbit
$\Cal O(\dash) \times \Cal O(\dash)$,
then we can argue, as in the proof of Theorem~A2 of [\cktrip],
that the proto-contact number \thetag{\triptwo}
correctly counts the number of simultaneous contacts.
In particular, we have the following result.
\par
\proclaim{Theorem~4}
Suppose that the curves $C$ and $D$ contain no lines,
and that the
general member of $\Cal X$ contains no line.
If $C$, $D$, and $\Cal X$ are in general position
with respect to the action of
$PGL(2) \times PGL(2)$
on $F(2) \times F(2)$,
then \thetag{\triptwo} counts the number of
simultaneous contacts
of order (3,3) between $C$ and $D$ and t
he members of $\Cal X$.
\endproclaim
\demo{Proof}
We will show that, for each nondense orbit
$\Cal O$ listed above, we have
$$
\dim\left((C(2) \times D(2)) \cap \Cal O \right)
	+ \dim \left(\Cal X_S(2,2) \cap (\Cal O \times S)\right)
	< \dim \Cal O .
\tag \dimineq
$$
Since $C$, $D$, and
$\Cal X$ are in general position with respect to the action,
transversality theory [\ktrans] then tells us that
the intersection
$$
(C(2) \times D(2)) \cap \sigma(\Cal X_S(2,2)) \cap \Cal O
$$
is empty.  Therefore the intersection of $\sigma(\Cal X_S(2,2))$
and $C(2) \times D(2)$ is transverse,
with all intersections
occurring in the dense orbit
$\Cal O(\dash) \times \Cal O(\dash)$.
By definition $\Cal X_S(2,2) = \Cal X(2) \fpos \Cal X(2)$,
where $\Cal X(2)$ is the closure of the graph of
a function defined
on a dense subset of $\Cal X$, and $C(2) \times D(2)$
is the closure of the graph of a function defined
on a dense subset of $C \times D$.
By a general position argument, all intersections between
$\sigma(\Cal X_S(2,2))$ and $C(2) \times D(2)$
are intersections between the graphs of these functions.
Hence each intersection point corresponds to a pair of points
$(c,d)$, in which $c$ is a nonsingular point of both $C$ and
some member $X_s$ of the family, and $d$ is a nonsingular
points of both $D$ and $X_s$.
The second-order data of $C$ and $X_s$ at $c$ are identical,
as are the second-order data of $D$ and $X_s$ at $d$.
Hence each intersection point is a simultaneous triple contact.
\par
To see that \thetag{\dimineq} holds for the nondense orbits,
note that the intersection of $C(2)$ or $D(2)$ with the divisor
at infinity is finite.  Thus
$$
\gather
\dim\left( (C(2) \times D(2))
	\cap(\Cal O(\dash) \times \Cal O(\infty))\right)
	\leq 1 \\
\dim\left( (C(2) \times D(2))
	\cap(\Cal O(\infty) \times \Cal O(\dash))\right)
	\leq 1 \\
\dim\left( (C(2) \times D(2))
	\cap(\Cal O(\infty) \times \Cal O(\infty))\right)
	\leq 0.
\endgather
$$
If $C$ and $D$  contain no lines, then $C(2)\cap \Cal O$
and $D(2)\cap \Cal O$ are each finite.  Hence
$$
\gather
\dim\left( (C(2) \times D(2))
	\cap(\Cal O(\dash) \times \Cal O(0))\right)
	\leq 1 \\
\dim\left( (C(2) \times D(2))
	\cap(\Cal O(0) \times \Cal O(\dash))\right)
	\leq 1 \\
\dim\left( (C(2) \times D(2))
	\cap(\Cal O(0) \times \Cal O(0))\right)
	\leq 0 \\
\dim\left( (C(2) \times D(2))
	\cap(\Cal O(0) \times \Cal O(\infty))\right)
	\leq 0 \\
\dim\left( (C(2) \times D(2))
	\cap(\Cal O(\infty) \times \Cal O(0))\right)
	\leq 0.
\endgather
$$
The dimension of $\Cal X_S(2,2)$ is 6, and, since
the general member of $\Cal X$ is reduced
(hence generically smooth), we obtain
$$
\gather
\dim\left( \Cal X_S(2,2)
	\cap(\Cal O(\dash) \times \Cal O(\infty) \times S)\right)
	\leq 5 \\
\dim\left( \Cal X_S(2,2)
	\cap(\Cal O(\infty) \times \Cal O(\dash) \times S)\right)
	\leq 5 \\
\dim\left( \Cal X_S(2,2)
	\cap(\Cal O(\infty) \times \Cal O(\infty) \times S)\right)
	\leq 5 \\
\dim\left( \Cal X_S(2,2))
	\cap(\Cal O(0) \times \Cal O(\infty) \times S)\right)
	\leq 5 \\
\dim\left( \Cal X_S(2,2))
	\cap(\Cal O(\infty) \times \Cal O(0) \times S)\right)
	\leq 5.
\endgather
$$
If the general member of $\Cal X$  contains no line,
then likewise
$$
\gather
\dim\left( \Cal X_S(2,2))
	\cap(\Cal O(\dash) \times \Cal O(0) \times S)\right)
	\leq 5 \\
\dim\left( \Cal X_S(2,2))
	\cap(\Cal O(0) \times \Cal O(\dash) \times S)\right)
	\leq 5 \\
\dim\left( \Cal X_S(2,2))
	\cap(\Cal O(0) \times \Cal O(0) \times S)\right)
	\leq 5. \\
\endgather
$$
Thus \thetag{\dimineq} is true in each of the eight cases.
\enddemo
\par
 Theorem~4 has a clear advantage over Theorem~3 in that
the hypotheses are readily verified.  Thus results like
Theorem~4 appear to be highly desirable.  To establish
such a result, however, one needs an understanding
of the $PGL(2)$-orbit structure of $F(n)$,
which becomes increasingly complicated as $n$ grows.
For example, there are 8 orbits on $F(3)$
and 21 orbits on $F(4)$.  (For a derivation of the first
number, see Theorem~2 of [\cktrip]. The second number
was first obtained
by Oberlin College student Dan Frankowski by {\it Mathematica}
calculations [\math] and later confirmed by us.)  Even worse,
since $PGL(2)$ has dimension~8, there are infinitely many orbits
on $F(n)$, none of them dense, when $n \geq 7$.
(In fact, the authors, along with Oberlin College students
Ian Robertson and Susan Sierra,
have found that there are infinitely many orbits on
$F(6)$.)
Hence theorems such as Theorem~4 cannot exist for
simultaneous contacts of arbitrary order.
This is why Theorem~3 is stated as it is, and appears
to be the best possible result of its type.
\par

\enddocument